%% file: MeshPaper.tex
\documentclass[12pt]{article}
\usepackage{amsfonts,amsmath,amssymb,amsthm}
\usepackage{a4,graphics,epsfig}

\usepackage[numbers,sort]{natbib}
\usepackage[table,xcdraw]{xcolor} 
\usepackage{float} 

\newcommand{\f}[1]{\mathbf{#1}}

\begin{document}
	
\begin{center}
{\LARGE \bf 
Isogeometric Parametrization Inspired by Large Elastic Deformation} \\[4mm]
Alexander Shamanskiy$^1$, Michael Helmut Gfrerer$^1$, Jochen Hinz$^2$ and Bernd Simeon$^1$ \\[1mm]
$^1$ TU Kaiserslautern, Dept. of Mathematics\\
$^2$ TU Delft, Dept. of Mathematics 
\end{center}

\begin{quote} \small
{\bf Abstract: }
The construction of volumetric parametrizations for computational domains is a key step in the pipeline of isogeometric analysis. Here, we investigate a solution to this problem based on the mesh deformation approach. The desired domain is modeled as a deformed configuration of an initial simple geometry. Assuming that the parametrization of the initial domain is bijective and that it is possible to find a locally invertible displacement field, the method yields a bijective parametrization of the target domain. We compute the displacement field by solving the equations of nonlinear elasticity with the neo-Hookean material law, and we show an efficient variation of the incremental loading algorithm tuned specifically to this application. In order to construct the initial domain, we simplify the target domain's boundary by means of an $L^2$-projection onto a coarse basis and then apply the Coons patch approach. The proposed methodology is not restricted to a single patch scenario but can be utilized to construct multi-patch parametrizations with naturally looking boundaries between neighboring patches. We illustrate its performance and compare the result to other established parametrization approaches on a range of two-dimensional and three-dimensional examples.

\end{quote}

\input{intro}

\input{deformation}

\input{method}

\input{initial}

\input{newton}

\input{examples}

\input{conclusion}

\bibliographystyle{unsrt}
\bibliography{refs}

\end{document}

%% file: intro.tex
\section{Introduction}
A common problem in isogeometric analysis (IGA) \cite{Hughes2005,Cottrell.2009} is generating a volumetric parametrization for the computational domain when only a description of its boundary is available. In this work, we investigate an approach to solving this problem which is based on mesh deformation. The parametrization for the target domain is acquired as a deformed configuration of a simple initial domain. The approach is related to a class of arbitrary Lagrangian-Eulerian methods in problems of fluid-structure interaction \cite{bazilevs2013,stein2003mesh,crosetto2011fluid} and to the interface tracking methods in free-surface flow problems \cite{zwicke2017boundary}. In the context of IGA, the approach has been applied in order to construct volumetric meshes consisting of a T-spline surface layer and a core of Lagrangian elements \cite{harmel2017volumetric}. Although only small deformations are considered, similar ideas are used to generate curvilinear meshes from piecewise linear triangulations in \cite{persson2009curved}.

We apply the mesh deformation approach to generate high-quality tensor product B-spline and NURBS parametrizations for complicated geometries.
It is done by first simplifying the target domain's boundary so that the Coons patch approach \cite{Piegl1997,farin1999discrete} can be applied to produce a bijective and uniform parametrization of the resulting simple geometry. The simplification can be conducted by means of projection in an $L^2$-sense onto a coarse basis; however, a number of ad hoc methods can be applied in every particular situation which makes the approach very flexible. Next, we deform the simplified geometry so that its boundary coincides with the target domain's boundary. We search for the unknown displacement field as a solution to the system of nonlinear elasticity equations with a prescribed boundary displacement. By using the logarithmic neo-Hookean material law, we exclude self-penetrations of the material and preserve the bijectivity of the initial parametrization. Moreover, we can partially preserve the uniformity of the initial parametrization by considering a nearly incompressible material. In order to efficiently solve the nonlinear elasticity equations, we employ a variation of the incremental loading algorithm. It numerically preserves bijectivity of the solution and can operate with an adaptive stepsize.   

The problem of generating tensor product B-spline and NURBS parametrizations has received a lot of attention since the introduction of IGA. Let us give a short overview of the state-of-the-art in the field. One of the simplest methods to construct a volumetric parametrization from a boundary description is the Coons patch. Although nothing guaranties that the resulting parametrization is bijective, the method is explicit, and its output can be used as a starting point for more sophisticated parametrization techniques. Nonlinear optimization \cite{Gravesen2014,xu2011parameterization,falini2015planar,xu2013constructing,pan2018low,ugalde2018injectivity} is a popular approach which allows to construct parametrizations optimal with respect to a chosen quality measure; the bijectivity is often enforced as an external constraint. Many other approaches seek to construct the parametrization as an inverse of a bijective mapping from the target domain to the parametric space. Among examples are the inverse of harmonic mappings \cite{nguyen2010parameterization,nguyen2012construction} and elliptic grid generation \cite{hinz2018elliptic}. For domains belonging to a certain class of geometries, specialized techniques have been developed. Examples are swept volumes \cite{aigner2009swept} and star-shaped domains \cite{Arioli2017}.

The rest of this paper is structured as follows. Section 2 gives a brief introduction into continuum mechanics and fixes necessary notation. Section 3 states the parametrization problem and then outlines the mesh deformation approach to its solution in broad brush-strokes. Section 4 deals with the construction of the initial domain, and numerical algorithms for computing the deformation are described in Section 5. A range of 2D and 3D examples is presented in Section 6 as well as a comparison of the results of the mesh deformation approach to other established parametrization techniques. Finally, Section 7 draws a conclusion and outlines further research directions. 

%% file: deformation.tex
\section{Nonlinear elasticity in a nutshell}
The following is a brief introduction into continuum mechanics $-$ based on \cite{wriggers2008nonlinear} $-$ where the focus lies on fixing a notation for the ingredients necessary for describing the mesh deformation approach in Section 3. 

Let $\Omega_0\subset\mathbb{R}^d$ be a reference configuration of a solid body undergoing a deformation $\pmb{\Phi}$. For each material point $\f{x}\in\Omega_0$, its position in the deformed configuration $\pmb{\Phi}(\f{x}) = \f{y}\in\Omega\subset\mathbb{R}^d$ can be expressed in terms of a displacement vector field $\f{u}:\Omega_0\to\mathbb{R}^d$ such that
\begin{equation}\label{eq:disp}
\f{y} = \f{x} + \f{u}(\f{x}).
\end{equation}  
Next, the deformation gradient $\f{F}:\Omega_0\to\mathbb{R}^{d\times d}$ is defined as
\begin{equation}
\f{F} = \nabla_\f{x}\pmb{\Phi} = \f{I} + \nabla_\f{x}\f{u}.
\end{equation}
Its determinant $J = \operatorname{det}\f{F}$ measures a relative volume change. Since self-penetration during deformation of the body is excluded, the mapping $\pmb{\Phi}$ must be bijective and the condition
\begin{equation}\label{eq:J}
J(\f{u}) > 0
\end{equation}
has to hold. In what follows we often use $J(\f{u})$ and $J(\pmb{\Phi})$ interchangeably.

Let the solid body be subject to volume forces $\f{g}$. From the conservation of linear momentum, it follows that the displacement $\f{u}$ fulfills the equations
\begin{equation}\label{eq:momentum}
-\operatorname{div}(\f{F}\f{S})(\f{u}) = \f{g} \text{ in }\Omega_0.
\end{equation} 
Here $\f{S}$ is the second Piola-Kirchhoff stress tensor which measures internal forces arising in the deformed solid body in response to the applied external load. Equations (\ref{eq:momentum}) are incomplete unless a relation between $\f{S}$ and $\f{u}$ $-$ called a material law $-$ is defined. In the present paper, we use
\begin{equation}\label{eq:neoHook}
\f{S}(\f{u}) = \lambda\ln J(\f{u})\f{C}(\f{u})^{-1} + \mu(\f{I}-\f{C}(\f{u})^{-1}),
\end{equation}
where the right Cauchy-Green tensor is defined as $\f{C} = \f{F}^T\f{F}$. The relation (\ref{eq:neoHook}) constitutes a particular choice of a nonlinear neo-Hookean material law. Note that due to the presence of $\ln J$, any displacement field $\f{u}$ satisfying the equations (\ref{eq:momentum}) with the material law (\ref{eq:neoHook}) grants a bijective deformation $\pmb{\Phi}$, i.e., (\ref{eq:J}) holds. The material law (\ref{eq:neoHook}) includes two constitutive parameters $-$ the so-called Lam\'e constants $\lambda$ and $\mu$ $-$ which can be computed from Young's modulus $E$ and Poisson's ratio $\nu$:
\begin{equation}
\lambda = \frac{\nu E}{(1+\nu)(1-2\nu)},\hspace{0.3cm} \mu = \frac{E}{2(1+\nu)}.
\end{equation} 
Finally, for the equations (\ref{eq:momentum}) to have a unique solution they have to be equipped with boundary conditions:
\begin{align}
\f{u} &= \f{u}_\mathcal{D} \text{ on } \partial\Omega_0^\mathcal{D},\\
\f{F}\f{S}\f{n} &= \f{f} \text{ on } \partial\Omega_0^\mathcal{N},
\end{align}
where $\partial\Omega_0^\mathcal{D}$ and $\partial\Omega_0^\mathcal{N}$ are the parts of the domain boundary $\partial\Omega_0$ with the prescribed displacement $\f{u}_\mathcal{D}$ and traction $\f{f}$; $\f{n}$ is the outer surface normal.

%% file: method.tex
\section{Mesh deformation approach}
Assume that for the target domain $\Omega\subset\mathbb{R}^d$ only a parametrization $\pmb{\partial}\f{G}(\pmb{\xi}):\partial[0,1]^d\to\partial\Omega$ of its boundary is available. The problem of domain parametrization is to construct a parametrization $\f{G}(\pmb{\xi}):[0,1]^d\to\Omega$ such that $\f{G}|_{\partial[0,1]^d} = \pmb{\partial}\f{G}$. Moreover, in order to be suitable for numerical simulations, the parametrization $\f{G}$ has to be bijective, i.e., the condition
\begin{equation}
J(\f{G}) = \operatorname{det}\nabla_{\pmb{\xi}}\f{G} > 0
\end{equation}
must hold. 

Assume further that the boundary parametrization $\pmb{\partial}\f{G}$ is given in terms of four compatible B-spline curves (for $d=2$) or six compatible surfaces (for $d=3$). By compatible we mean that the oppositely lying parts of the boundary have the same B-spline basis. In this case, the tensor product basis $\{B_i(\pmb{\xi})\}$ of the unknown parametrization $\f{G}$ is defined, and $\f{G}$ has the structure
\begin{equation}\label{eq:basis}
\f{G}(\pmb{\xi}) = \sum_{i=1}^{n}\f{c}_iB_i(\pmb{\xi}),
\end{equation}
where $\{\f{c}_i\}_{i=1}^n$ are the control points. Since the boundary control points $\{\f{c}_i\}_\mathcal{B}$ follow from $\pmb{\partial}\f{G}$, the problem boils down to allocation of the unknown interior control points $\{\f{c}_i\}_\mathcal{I}$.

We apply the mesh deformation approach to solve the stated parametrization problem. The idea is to start by choosing a simple initial domain $\Omega_0$ with a known parametrization $\f{G}_0:[0,1]^d\to\Omega_0$. We assume that the parametrization $\f{G}_0$ uses the same tensor product basis as $\f{G}$ and thus has the following form:
\begin{equation}
\f{G}_0(\pmb{\xi}) = \sum_{i=1}^n\f{c}_i^0 B_i(\pmb{\xi}).
\end{equation}
Next, we search for a deformation $\pmb{\Phi}:\Omega_0\to\Omega$ such that 
\begin{equation}\label{eq:bdryfit}
\pmb{\Phi}(\pmb{\partial}\f{G}_0(\pmb{\xi})) = \pmb{\partial\f{G}}(\pmb{\xi}) \text{ for } \forall\pmb{\xi}\in\partial[0,1]^d.
\end{equation}
The deformation $\pmb{\Phi}$ is characterized by an unknown displacement field $\f{u}:\Omega_0\to\mathbb{R}^d$. Following the isogeometric approach, we can introduce the discretization
\begin{equation}\label{eq:uh}
\f{u}_h(\f{x}) = \sum_{i=1}^n\f{d}_i B_i(\f{G}_0^{-1}(\f{x})),
\end{equation}
where the boundary degrees of freedom $\{\f{d}_i\}_\mathcal{B}$ are given by (\ref{eq:bdryfit}) as
\begin{equation}\label{eq:DBC2}
\{\f{d}_i\}_\mathcal{B} =\{\f{c}_i-\f{c}_i^0\}_\mathcal{B}.
\end{equation}
Once the interior degrees of freedom $\{\f{d}_i\}_\mathcal{I}$ are found, 
the parametrization $\f{G}$ can be constructed as a composition of the parametrization $\f{G}_0$ and the deformation $\pmb{\Phi}$:
\begin{equation}\label{eq:super}
\f{G} = \pmb{\Phi}\circ\f{G}_0
\end{equation}
or
\begin{equation}\label{eq:G}
\f{G}(\pmb{\xi}) =  \f{G}_0(\pmb{\xi}) + \f{u}_h(\f{G}_0(\pmb{\xi}))) = \sum_{i=1}^{n}(\f{c}_i^0+\f{d}_i)B_i(\pmb{\xi}).
\end{equation}
Observe that if the initial parametrization $\f{G}_0$ and the deformation $\pmb{\Phi}$ are bijective, i.e.,
\begin{equation}
J(\f{G}_0)>0\; \text{ and }\; J(\pmb{\Phi}) > 0
\end{equation}
hold, then the resulting parametrization $\f{G}$ is bijective as well: 
\begin{equation}\label{eq:superJ}
J(\f{G}) = \operatorname{det}\nabla_{\pmb{\xi}}\f{G} = \operatorname{det}\nabla_\f{x}\pmb{\Phi} \operatorname{det}\nabla_{\pmb{\xi}}\f{G}_0  = J(\pmb{\Phi})J(\f{G}_0) >0. 
\end{equation}

The choice of the initial domain $\Omega_0$ is discussed in Section 4. To compute the deformation $\pmb{\Phi}$, we use the equations of nonlinear elasticity introduced in Section 2. The initial domain $\Omega_0$ serves as a reference configuration, and the unknown displacement field $\f{u}$ is found as a solution to the following system of equations:
\begin{align}\label{eq:nonlin}
-\operatorname{div}(\f{F}\f{S})(\f{u}) &= \f{0} \text{ in } \Omega_0,\\
\f{u} &= \f{u}_\mathcal{D} \text{ on } \partial\Omega_0,\label{eq:DBC}
\end{align}
where $\f{u}_\mathcal{D}$ is the prescribed boundary displacement defined by the boundary degrees of freedom (\ref{eq:DBC2}):
\begin{equation}
\f{u}_\mathcal{D}(\f{x}) = \sum_{i\in\mathcal{B}}\f{d}_i B_i(\f{G}_0^{-1}(\f{x})).
\end{equation}

As for the material parameters, the choice of Young's modulus $E$ does not affect the solution of the system (\ref{eq:nonlin}-\ref{eq:DBC}) since equation (\ref{eq:nonlin}) has a zero right-hand side and the Dirichlet boundary condition (\ref{eq:DBC}) is prescribed over the entire boundary of the domain $\Omega_0$. On the other hand, Poisson's ratio $\nu$ is of great importance since it determines the resistance of the material to volumetric changes. A material with high Poisson's ratio will resist self-penetration and will thus contribute to the preservation of bijectivity. When $\nu$ approaches 0.5, the material becomes nearly incompressible. In practice, we use values between 0.45 and 0.49 since values higher than 0.49 would lead to a numerically unstable system unless a special formulation for incompressible behavior is used. A truly incompressible material, though,  does not allow for any volumetric changes and would require the domains $\Omega_0$ and $\Omega$ to be of the same volume $-$ a condition which is hard to satisfy in practice.

%% file: initial.tex
\section{Initial domain}
The choice of the initial domain $\Omega_0$ is a rather empirical step which directly affects the quality of the resulting parametrization $\f{G}$ for the target domain $\Omega$. Ideally, $\Omega_0$ should be simple enough $-$ so that it is possible to parametrize it using the Coons patch approach or another explicit method $-$ and yet geometrically close enough to $\Omega$  $-$ so that the complexity of computing a bijective deformation $\pmb{\Phi}$ does not eclipse the complexity of the original parametrization problem for $\Omega$. A combination of these two requirements suggests that the boundary of $\Omega_0$ has to be a simplification of the target domain's boundary $\partial\Omega$. In what follows, we describe a basic simplification procedure which allows to generate a range of different initial domains. 

\subsection{Boundary simplification}
We propose a simplification technique which is based on a projection onto a coarse B-spline basis in the $L^2$-sense.
The idea is that only geometrically simple shapes can lie in the span of such a basis. The boundary of the target domain can be simplified as a whole or in parts. Due to a tensor product structure of the parametrization $\f{G}$, it is convenient to simplify each side of the domain separately; one should only make sure that the simplified sides fit together at the interfaces. After projection, the simplified boundary is re-expressed in terms of the original basis.  

Let $\Gamma=\pmb{\partial}\f{G}|_\Pi$ be a parametrization of the part of the target domain boundary $\partial\Omega$ where $\Pi\subset\partial[0,1]^d$. Additionally, let $\mathcal{P}$ denote a set of indices corresponding to basis functions $B_i(\xi)$ (\ref{eq:basis}) which are not zero on $\Pi$. Then 
\begin{equation}
\Gamma(\pmb{\xi}) = \sum_{i\in\mathcal{P}}\f{c}_iB_i(\pmb{\xi}).
\end{equation}
In order to construct a simplification of $\Gamma$, we introduce a coarse basis $\{b_i(\pmb{\xi})\}_{i=1}^m$ where $m\ll|\mathcal{P}|$. By projecting $\Gamma$ in the $L^2$-sense onto the basis $\{b_i(\pmb{\xi})\}_{i=1}^m$, we acquire a primary simplification $\gamma$,
\begin{equation}\label{eq:prsimple}
\gamma(\xi) = \sum_{i=1}^{m}\f{x}_ib_i(\pmb{\xi}).
\end{equation}
The control points $\{\f{x}_i\}_{i=1}^m$ are found by solving the linear system
\begin{align}\label{eq:L2EQ}
\Big((b_i,b_j)_\Pi\Big)\Big(\f{x}_i^T\Big)&=\Big((b_i,\Gamma)_\Pi^T\Big),\\
\gamma|_{\partial\Pi} &= \Gamma,\label{eq:L2DBC}
\end{align}
where the inner product $(A,B)_\Pi$ is defined as $\int_{\Pi}A(\pmb{\xi})B(\pmb{\xi})d\pmb{\xi}.$ The boundary condition (\ref{eq:L2DBC}) ensures that the simplifications of different parts of $\partial\Omega$ fit together at $\partial\Pi$. 

The primary simplification $\gamma$ can be re-expressed in terms of the original basis $\{B_i(\pmb{\xi})\}_\mathcal{P}$ in two ways: either by applying h- and p-refinement $-$ also known as knot insertion and degree elevation $-$ or by projecting $\gamma$ onto $\{B_i(\pmb{\xi})\}_{i\in\mathcal{P}}$ in a manner analogous to (\ref{eq:L2EQ}-\ref{eq:L2DBC}). The latter slightly changes the shape of $\gamma$ which is insignificant since we have freedom in choosing the initial domain $\Omega_0$. The result is a simplification $\Gamma_0$:
\begin{equation}
\Gamma_0(\pmb{\xi}) = \sum_{i\in\mathcal{P}}\f{c}_i^0B_i(\pmb{\xi}).
\end{equation}

\begin{figure}[H]
	\centering
	\includegraphics[height = 3cm]{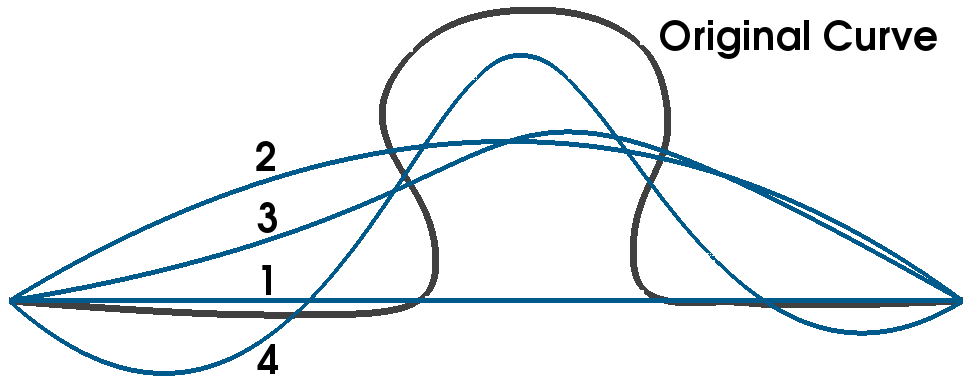}
	\caption{Dependence of the simplified boundary on the polynomial degree of the coarse basis. The coarsest B-spline bases of degree 1, 2, 3, and 4 are used.}\label{fig:simplify1}
\end{figure}

\begin{figure}[H]
	\centering
	\includegraphics[height = 3cm]{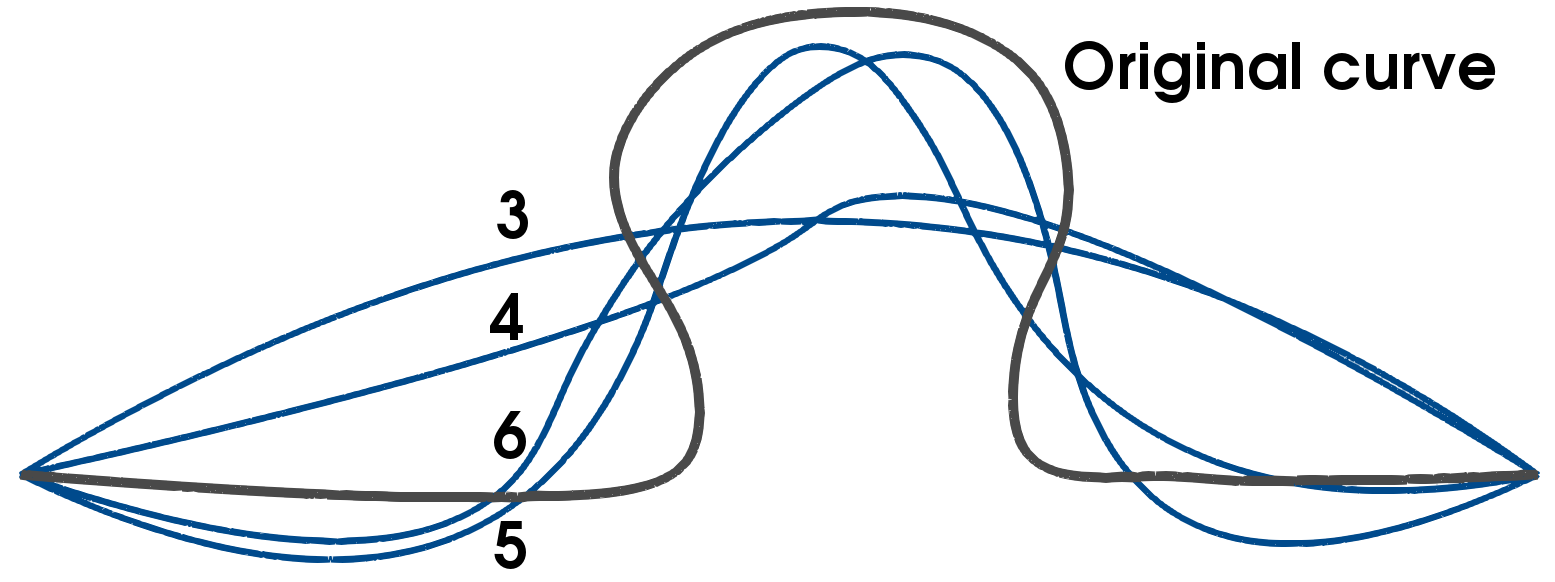}
	\caption{Dependence of the simplified boundary on the number of basis functions. Quadratic B-spline bases with 3, 4, 5, and 6 elements are used.}\label{fig:simplify2}
\end{figure}

The actual shape of $\Gamma_0$ depends on the choice of the coarse basis $\{b_i(\pmb{\xi})\}_{i=1}^m$. Figure~\ref{fig:simplify1} shows an example where the coarsest B-spline bases of degree 1, 2, 3, and 4 are used for projection. Observe the rapid growth of complexity of the resulting simplified geometry as the polynomial degree increases. The same effect is observed with the growth of the number of basis functions $m$ with a fixed polynomial degree, see Fig.~\ref{fig:simplify2}.  

One of the advantages of the proposed simplification technique is a partial preservation of the parametrization speed which means that the images $\Gamma(\pmb{\xi})$ and $\Gamma_0(\pmb{\xi})$ of the same parametric point $\pmb{\xi}$ are close to each other. This reduces the prescribed boundary displacement in (\ref{eq:DBC}), makes it easier to compute the deformation $\pmb{\Phi}$, and thus increases the quality of the resulting parametrization $\f{G}$. 

\subsection{Coons patch}
The procedure described above is applied to the entire boundary $\partial\Omega$. The result is the parametrization of the initial domain's boundary $\pmb{\partial}\f{G}_0:\partial[0,1]^d\to\partial\Omega_0$:
\begin{equation}
\pmb{\partial}\f{G}_0(\pmb{\xi}) = \sum_{i\in\mathcal{B}}\f{c}_i^0B(\pmb{\xi}).
\end{equation}
Our intention is to parametrize it using the Coons patch approach. In a two-dimensional case, the Coons patch defines $\f{G}_0$ as a bilinear blending of four parametric curves $\pmb{\partial}\f{G}_0(0,\xi_2)$, $\pmb{\partial}\f{G}_0(1,\xi_2)$, $\pmb{\partial}\f{G}_0(\xi_1,0)$ and $\pmb{\partial}\f{G}_0(\xi_1,1)$,
\begin{align}
\f{G}_0(\xi_1,\xi_2) = &(1-\xi_1)\pmb{\partial}\f{G}_0(0,\xi_2) + \xi_1\pmb{\partial}\f{G}_0(1,\xi_2)  \nonumber\\
\label{eq:coons}+&(1-\xi_2)\pmb{\partial}\f{G}_0(\xi_1,0) + \xi_2\pmb{\partial}\f{G}_0(\xi_1,1) \\ 
-&\begin{bmatrix}
1-\xi_1  & \xi_1  
\end{bmatrix}
\begin{bmatrix}
\pmb{\partial}\f{G}_0(0,0) & \pmb{\partial}\f{G}_0(0,1)  \\
\pmb{\partial}\f{G}_0(1,0) & \pmb{\partial}\f{G}_0(1,1) 
\end{bmatrix}
\begin{bmatrix}
1-\xi_2  \\
\xi_2  
\end{bmatrix},\nonumber
\end{align}
where $\xi_1$ and $\xi_2$ are parametric coordinates. The provided definition can be straightforwardly generalized to a three-dimensional case. 

\begin{figure}[H]
	\centering
	\includegraphics[height = 3.95cm]{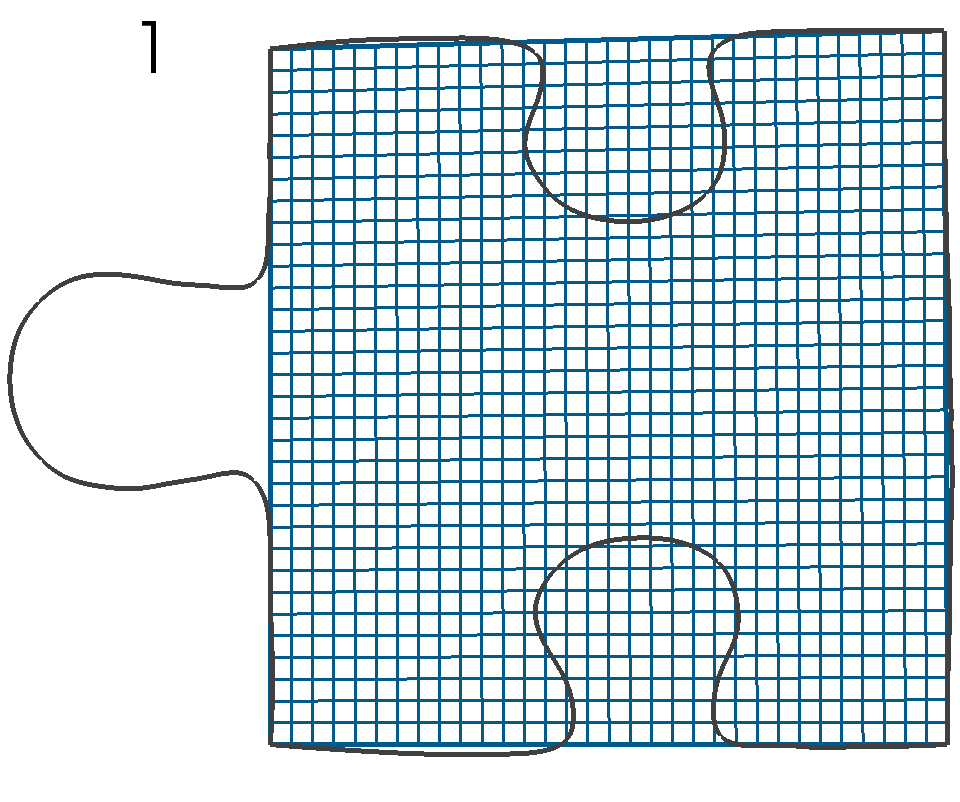}
	\includegraphics[height = 3.95cm]{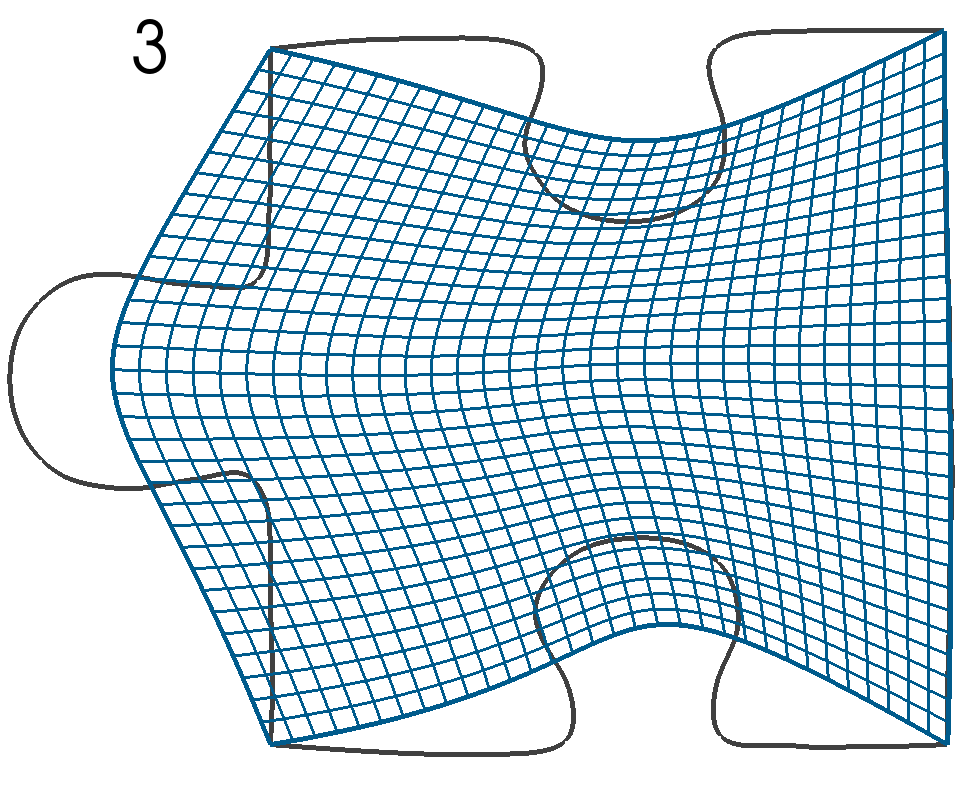}
	\includegraphics[height = 3.95cm]{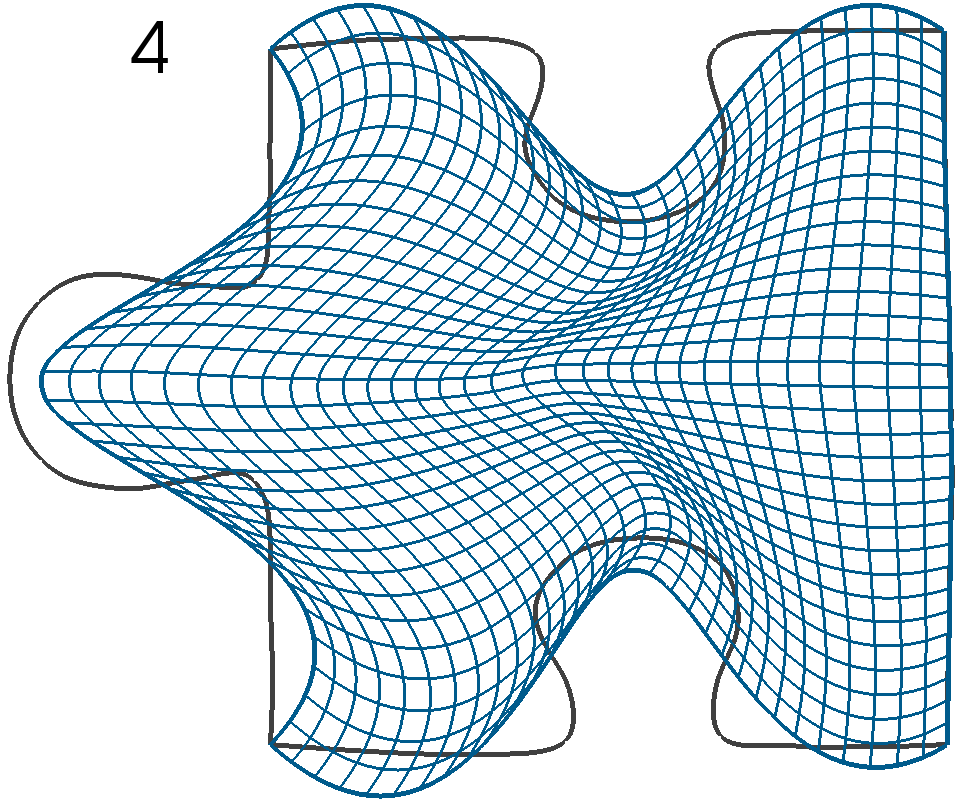}
	\caption{Coons patch applied to different simplifications of the puzzle piece geometry. The coarsest B-spline bases are of degree 1, 3, and 4 are used for simplification.}\label{fig:coons}
\end{figure}
Note that nothing guarantees that the resulting parametrization $\f{G}_0$ is bijective. We assume, however, that the boundary of the initial domain $\Omega_0$ acquired at the simplification step is simple enough for the Coons patch to succeed. This assumption puts a restriction on how fine the coarse basis (\ref{eq:prsimple}) can be. If the coarsest linear basis is used $-$ which is equivalent to substituting the original domain by a quad spanned on its corners $-$, then the Coons patch approach always produces a uniform bijective parametrization if the quad is convex. On the other hand, the quad-simplification is often too simple. There may exist a different initial domain which can also be parametrized with the Coons patch but which is geometrically closer to the target domain, see Fig.~\ref{fig:coons}. In the case of the depicted puzzle piece example, an optimal simplification is acquired by using the coarsest cubic basis.

%% file: newton.tex
\section{Deformation}
\subsection{Incremental Newton's method}
The nonlinear system (\ref{eq:nonlin}-\ref{eq:DBC}) is usually solved using Newton's method. However, if the prescribed boundary displacement $\f{u}_\mathcal{D}$ (\ref{eq:DBC}) is large, it can be difficult to find a bijective initial guess satisfying (\ref{eq:DBC}) from where Newton's method could converge to the system's solution $\f{u}$. In this case, the incremental loading can be applied, i.e., the problem (\ref{eq:nonlin}-\ref{eq:DBC}) is replaced by a sequence of problems for each loading step $i=1,\dots,N$:
\begin{align}\label{eq:nonlininc}
-\operatorname{div}(\f{F}\f{S})(\f{u}^i) &= \f{0} \text{ in } \Omega_0,\\
\f{u}^i &= \frac{i}{N}\f{u}_\mathcal{D} \text{ on } \partial\Omega_0.\label{eq:DBCinc}
\end{align}
Each incremental displacement $\f{u}^i$ provides an initial guess for Newton's method at the next loading step. 

To formulate the algorithm, we need a weak form of equations (\ref{eq:nonlininc}-\ref{eq:DBCinc}). Let $\mathcal{V}_i = \{\f{v}\in H^1(\Omega_0)^d \;|\; \f{v} = \frac{i}{N}\f{u}_\mathcal{D} \text{ on }\partial\Omega_0\}$ be a set of trial solution spaces for different loading steps $i=1,\dots,N$ and let the weighting function space $\mathcal{V}_0$ be defined as $\{\f{v}\in H^1(\Omega_0)^d \;|\; \f{v} = \f{0}\text{ on }\partial\Omega_0\}$. Then the weak form of equations (\ref{eq:nonlininc}-\ref{eq:DBCinc}) is
\begin{align}\nonumber
&Find\text{ } \f{u}^i\in\mathcal{V}_i \text{ }such\text{ } that\\
&P(\f{u}^i,\f{v}) = \int\displaylimits_{\Omega_0}\f{S}(\f{u}^i):\delta\f{E}(\f{u}^i)[\f{v}]d\f{x} = 0, \hspace{0.2cm} \forall \f{v}\in\mathcal{V}_0,\label{eq:weak}
\end{align}
where $\delta\f{E}(\f{u}^*)[\f{v}] = \frac{1}{2}\big(\f{F}(\f{u}^*)^T\nabla_\f{x}\f{v} + \nabla_\f{x}\f{v}^T\f{F}(\f{u}^*)\big)$ is the variation of the Green-Lagrange strain tensor. Equations (\ref{eq:weak}) are nonlinear; in order to apply Newton's method they have to be linearized. The Taylor expansion at $P(\f{u}^*,\f{v})$ with the displacement increment $\Delta\f{u}$ yields 
\begin{equation}
P(\f{u}^*+\Delta\f{u},\f{v}) = P(\f{u}^*,\f{v}) + DP(\f{u}^*,\f{v})\cdot\Delta\f{u} + o(||\Delta\f{u}||),
\end{equation}
where the directional derivative $DP(\f{u}^*,\f{v})\cdot\Delta\f{u}$ is given by
\begin{equation}\label{eq:dir}
DP(\f{u}^*,\f{v})\cdot\Delta\f{u} = \int\displaylimits_{\Omega_0}\Big(\nabla_\f{x}\Delta\f{u}\,\f{S}(\f{u}^*):\nabla_\f{x}\f{v} +\mathbb{C}(\f{u}^*)\delta\f{E}(\f{u}^*)[\Delta\f{u}]:\delta\f{E}(\f{u}^*)[\f{v}]\Big)d\f{x}.
\end{equation}
Here $\mathbb{C}=2\frac{d\f{S}}{d\f{C}}$ is the forth order elasticity tensor whose components, in case of the neo-Hookean material law, are given by
\begin{equation}
\mathbb{C}_{abcd} = \lambda\f{C}_{ab}^{-1}\f{C}_{cd}^{-1}+(\mu-\lambda\ln J)\big(\f{C}_{ac}^{-1}\f{C}_{bd}^{-1}+\f{C}_{ad}^{-1}\f{C}_{bc}^{-1}\big).
\end{equation}

Having defined all the necessary tools, we can now formulate incremental Newton's method. Let $\f{u}_s^i$ be a displacement field at the $s$-th iteration of Newton's method at the $i$-th loading step. The method involves two operation types. Type-A is an update $\f{u}_s^i\in\mathcal{V}_i\to\f{u}_{s+1}^i\in\mathcal{V}_i$ within the $i$-th loading step. An increment $\Delta\f{u}_{s}^i\in\mathcal{V}_0$ such that $\f{u}_{s+1}^i=\f{u}_{s}^i+\Delta\f{u}_{s}^i$ is found as the solution to the following weak problem: 
\begin{align}\nonumber
&Find\text{ }\Delta\f{u}_{s}^i\in\mathcal{V}_0\text{ }such\text{ }that\\
\label{eq:incweak2}&DP(\f{u}_{s}^i,\f{v})\cdot\Delta\f{u}_{s}^{i} = - P(\f{u}_{s}^i,\f{v}), \hspace{0.2cm} \forall \f{v}\in\mathcal{V}_0.
\end{align}
We repeat this step untill the convergence criterion
\begin{equation}
\frac{||\Delta\f{u}_s^i||_{L^2}}{||\f{u}_s^i||_{L^2}} < \varepsilon
\end{equation}
is met. The last approximate solution defines the incremental displacement $\f{u}^i$ at the $i$-th loading step.

Type-B is an update $\f{u}^{i-1}\in\mathcal{V}_{i-1}\to\f{u}_{1}^{i}\in\mathcal{V}_{i}$ between loading steps. We search for an increment $\Delta\f{u}^{i}\in\mathcal{V}_1$ such that $\f{u}_{1}^{i}=\f{u}^{i-1}+\Delta\f{u}^{i}$ as a solution to the weak problem
\begin{align}\nonumber
&Find\text{ }\Delta\f{u}^{i}\in\mathcal{V}_1\text{ }such\text{ }that\\
\label{eq:incweak4}&DP(\f{u}^{i-1},\f{v})\cdot\Delta\f{u}^{i} = - P(\f{u}^{i-1},\f{v}), \hspace{0.2cm} \forall \f{v}\in\mathcal{V}_0.
\end{align}
We say that the increment $\Delta\f{u}_i$ has stepsize $h_i = 1/N$ meaning that $\Delta\f{u}_i$ advances the displacement at the boundary $\partial\Omega_0$ by $1/N$-th of $\f{u}_\mathcal{D}$ (\ref{eq:DBC}).

The method is initialized with an initial displacement $\f{u}^0 = \f{0}$. The incremental displacement $\f{u}^N$ at the $N$-th loading step is accepted as the approximate solution $\f{u}$ to the original system (\ref{eq:nonlin}-\ref{eq:DBC}).

\subsection{Bijectivity and adaptivity}
Although the material law (\ref{eq:neoHook}) guarantees that the solution $\f{u}$ to the system (\ref{eq:nonlin}-\ref{eq:DBC}) is bijective, special care is required to achieve this property when solving the system using Newton's method. Note that the directional derivative $DP(\f{u},\f{v})\cdot\Delta\f{u}$ (\ref{eq:dir}) can only be evaluated at a bijective displacement $\f{u}^*$. However, both type-A (\ref{eq:incweak2}) and type-B (\ref{eq:incweak4}) updates can produce an increment $\Delta\f{u}$ such that $\f{u}^*+\Delta\f{u}$ is not bijective. The problem can be overcome by adaptively scaling the increment $\Delta\f{u}$. If $J(\f{u}^*) > 0$, there always exists a scaling coefficient $t\in[0,1]$ such that 
\begin{equation}\label{eq:scaling}
J(\f{u}^*+t\Delta\f{u}) > 0.
\end{equation}
In practice, we determine the scaling coefficient $t$ by consecutively testing values $t^k = 1/2^k$ until (\ref{eq:scaling}) is satisfied.

The implementation of adaptivity differs slightly for the type-A and type-B updates. For a type-A update $\f{u}_{s+1}^i = \f{u}_s^i+\Delta\f{u}_{s+1}^i$ $-$ where the increment $\Delta\f{u}_s^i$ is determined solely by (\ref{eq:incweak2}) $-$ the scaled increment $t\Delta\f{u}_s^i$ simply redefines $\f{u}_{s+1}^i$ as $\f{u}_s^i + t\Delta\f{u}_s^i$, and the method proceeds to the next iteration of Newton's method. 

For a type-B update $\f{u}_1^i = \f{u}^{i-1}+\Delta\f{u}^i$, the stepsize $h_i$ of the increment $\Delta\f{u}^i$ is predefined by the number of loading steps $N$. Scaling the increment $\Delta\f{u}^i$ changes the stepsize to $t\cdot h_i$ and $-$ since all updates of the boundary displacement have to add up to $\f{u}_\mathcal{D}$ $-$ requires changing the stepsizes of the subsequent type-B updates. One way to do it is to proceed with the $1/N$ stepsize, scaling it if necessary to fulfill (\ref{eq:scaling}). The final stepsize
\begin{equation}
h_{N^*} = 1- \sum_{i=1}^{N^*-1}h_i,
\end{equation}
makes sure that all stepsizes add up to 1.

Another possibility is to apply the so-called greedy stepsize strategy where incremental Newton's method begins with a stepsize of the first type-B update $h_1=1$. If the resulting displacement $\f{u}_1$ is not bijective, $h_1$ is iteratively halved until (\ref{eq:scaling}) is satisfied. The method proceeds with the stepsize
\begin{equation}
h_i = 1 - \sum_{j=1}^{i-1}h_j
\end{equation}
for loading steps $i\geqslant2$ which is also iteratively halved if necessary. We have to report that the greedy stepsize often results in stalling of the adaptive algorithm, i.e., the iterative halving produces too small stepsizes. The effect occurs much less often with the first adaptive strategy. This behavior deserves further investigation.

A nonadaptive solution to preserve bijectivity during type-B updates is to increase the number of loading steps $N$ and to restart the method.

Lastly, we remark on ways to test the bijectivity condition (\ref{eq:J}). A solution which takes into account the B-spline nature of the discretization $\f{u}_h$ (\ref{eq:uh}) is to express the Jacobian determinant $J(\f{u}_h)$ as a B-spline function \cite{Gravesen2014}. If all control coefficients in a B-spline expansion of $J(\f{u}_h)$ are positive, then the displacement $\f{u}_h$ is bijective. Unfortunately, this condition is only a sufficient but not a necessary one; this may lead to a lot of false detections of bijectivity violation. In practice, we resort to a much less elegant solution of sampling the Jacobian determinant at the Gaussian quadrature points associated with the discretization $\f{u}_h$. 
\subsection{Diagonal incremental loading}
Notice that incremental Newton's method is computationally expensive. If $S$ is the average number of iterations which Newton's method takes to converge at each loading step, then the method requires $O(NS)$ iterations to compute $\f{u}$. This is justified for applications where the deformation history is important; in our case, however, only the final displacement field $\f{u}^N$ is of interest.

In what follows, we propose a variation of incremental Newton's method which requires only $O(N+S)$ iterations. It begins with $N$ type-B updates $\f{u}_{inc}^{i-1}\in\mathcal{V}_{i-1}\to \f{u}_{inc}^{i}\in\mathcal{V}_{i}$ between loading steps. Similar to (\ref{eq:incweak4}), an increment $\Delta\f{u}_{inc}^{i}\in\mathcal{V}_1$ such that $\f{u}_{inc}^{i}=\f{u}_{inc}^{i-1}+\Delta\f{u}_{inc}^{i}$ is found as a solution to the weak problem
\begin{align}\nonumber
&Find\text{ }\Delta\f{u}_{inc}^{i}\in\mathcal{V}_1\text{ }such\text{ }that\\
\label{eq:incweak6}&DP(\f{u}_{inc}^{i-1},\f{v})\cdot\Delta\f{u}_{inc}^{i} = - P(\f{u}_{inc}^{i-1},\f{v}), \hspace{0.2cm} \forall \f{v}\in\mathcal{V}_0.
\end{align}
Once again, the method is initialized with $\f{u}_{inc}^0=\f{0}$. After $N$ steps, the displacement $\f{u}_{inc}^N\in\mathcal{V}_N$ is acquired; the described above adaptive algorithms can be applied. From here, the method proceeds with type-A iterations (\ref{eq:incweak2}) till it converges to $\f{u}^N=\f{u}$.
\begin{figure}[H]
	\centering
	\includegraphics[height = 3.6cm]{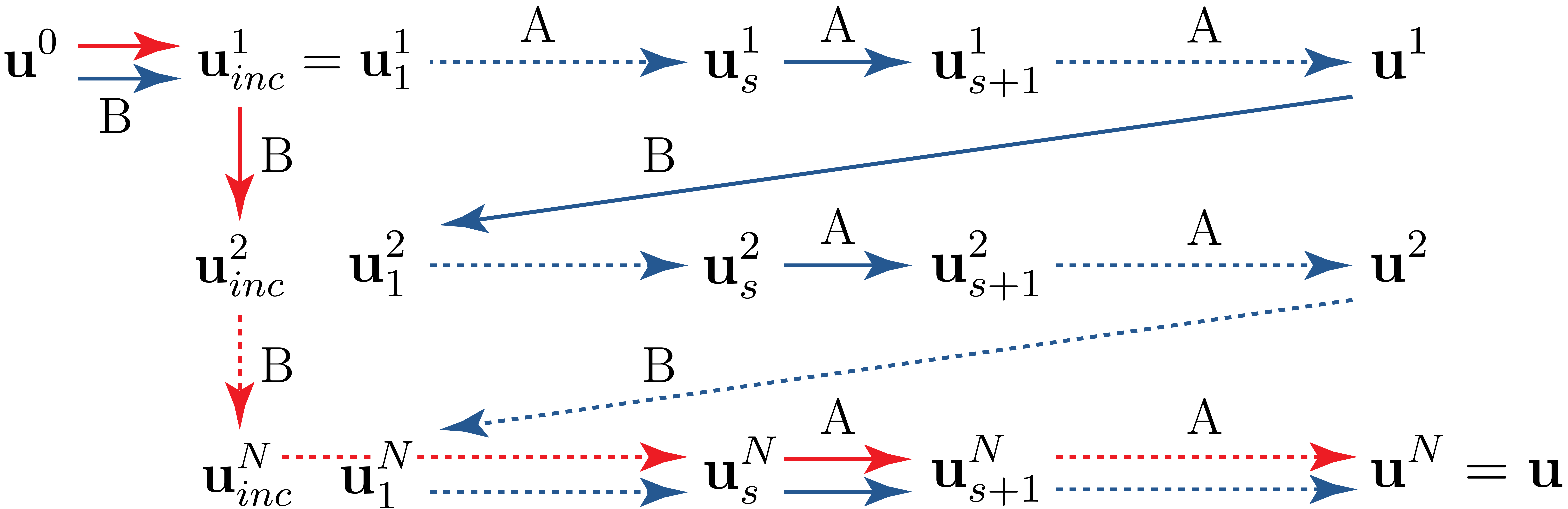}
	\caption{Incremental Newton's method (blue) and the proposed variation with diagonal incremental loading (red).}\label{fig:newton}
\end{figure}

The difference between the incremental Newton's method and the proposed variation is illustrated schematically in Figure~\ref{fig:newton}. We refer to the first phase of the algorithm as the diagonal incremental loading since the updates $\f{u}_{inc}^{i-1}\to \f{u}_{inc}^{i}$ advance the solution through both the iterations of Newton's method and the loading steps. In fact, our numerical experiments suggest that $\f{u}_{inc}^N$ converges quadratically to $\f{u}$ as $N\to\infty$. Because of that, $\f{u}_{inc}^N$ can be used as a stand-alone approximate solution to the system (\ref{eq:nonlin}-\ref{eq:DBC}) if the number of loading steps $N$ is big enough. In this case, the use of an adaptive loading stepsize is unnecessary. 

\subsection{Diagonal incremental loading with linear elasticity}
The incremental displacements $\f{u}^i$ define a sequence of intermediate domains $\Omega_i = \{\f{x}+\f{u}^i(\f{x})\; | \;\f{x}\in\Omega_0\}$, $i=1,\dots,N$. If the number of loading steps $N$ is big enough, one could try to construct a similar sequence $\Omega_i^{lin} = \{\f{x}+\Delta\f{u}^i_{lin}(\f{x})\; | \;\f{x}\in\Omega_{i-1}^{lin}\}$, $\Omega_0^{lin}=\Omega_0$ recursively, where at the $i$-th step a displacement increment $\Delta\f{u}_{lin}^i$ is found as a solution to the system of linear elasticity equations
\begin{align}\label{eq:linel}
-\operatorname{div}\pmb{\sigma}(\Delta\f{u}^i_{lin}) &= \f{0}\text{ in } \Omega_{i-1}^{lin},\\
\label{eq:linelDBC}\Delta\f{u}^i_{lin} &= \frac{\f{u}_\mathcal{D}}{N}\text{ on }\partial\Omega_{i-1}^{lin}.
\end{align}
Here $\pmb{\sigma}$ is the Cauchy stress tensor which is related to the linear strain tensor $\pmb{\varepsilon}(\f{u}^*)=\frac{1}{2}\big(\nabla_\f{x}\f{u}^*+(\nabla_\f{x}\f{u}^*)^T\big)$ via Hooke's law:
\begin{equation}
\pmb{\sigma}(\f{u}^*) = \lambda\operatorname{tr}(\pmb{\varepsilon}(\f{u}^*))\f{I} + 2\mu\pmb{\varepsilon}(\f{u}^*). 
\end{equation}

Note that, unlike (\ref{eq:nonlininc}-\ref{eq:DBCinc}), the equations (\ref{eq:linel}-\ref{eq:linelDBC}) are formulated in the intermediate configurations $\Omega_i^{lin}$, not in $\Omega_0$. For each $\Omega_i^{lin}$ we define the trial solution space $\mathcal{V}_i^{lin} = \{\f{v}\in H^1(\Omega_{i-1}^{lin})^d \;|\; \f{v} = \frac{\f{u}_\mathcal{D}}{N}\text{ on }\partial\Omega_{i-1}\}$ and the weighting function space $\mathcal{V}_{i,0}^{lin} = \{\f{v}\in H^1(\Omega_{i-1}^{lin})^d \;|\; \f{v} = \f{0} \text{ on }\partial\Omega_{i-1}\}$. Then the weak form of equations (\ref{eq:linel}-\ref{eq:linelDBC}) is
\begin{align}\nonumber
&Find\text{ }\Delta\f{u}^i_{lin}\in\mathcal{V}_i^{lin}\text{ }such\text{ }that\\
\label{eq:incweak8}&L(\Delta\f{u}^i_{lin},\f{v})=\int\displaylimits_{\Omega_{i-1}}\pmb{\sigma}(\Delta\f{u}^i_{lin}):\pmb{\varepsilon}(\f{v})d\f{x} = 0, \hspace{0.2cm} \forall \f{v}\in\mathcal{V}_{i,0}^{lin}.
\end{align}

It is important to notice that the weak problem (\ref{eq:incweak8}) is not equivalent to (\ref{eq:incweak6}). In fact, the bilinear form $L(\Delta\f{u},\f{v})$ in (\ref{eq:incweak8}) is the result of evaluating the directional derivative $DP(\f{u}^*,\f{v})\cdot\Delta\f{u}$ (\ref{eq:dir}) at $\f{u}^*=\f{0}$. Thus, the described procedure $-$ which we refer to as the linear diagonal incremental loading as opposed to the described above (nonlinear) diagonal incremental loading $-$ is similar to modified Newton's method in \cite{wriggers2008nonlinear} where the derivative evaluated at the first iteration is used to compute updates at all consecutive iterations.

We define the linear incremental displacements $\f{u}_{lin}^i$ as a sum of the preceding increments $\Delta\f{u}_{lin}^i$:
\begin{equation}
\f{u}_{lin}^i=\sum_{j=1}^i\Delta\f{u}_{lin}^j.
\end{equation}
In our experience, as the number of loading steps $N$ grows, $\f{u}_{lin}^N$ converges linearly to a displacement $\f{u}_{lin}$ which, although not equal, is close to the solution $\f{u}$ of the system (\ref{eq:nonlin}-\ref{eq:DBC}). Even more importantly, as we demonstrate in Section 6, the limit displacement $\f{u}_{lin}$ seems to be bijective; the described above adaptive algorithms can be applied to ensure bijectivity for small $N$.

Much like diagonal incremental loading with nonlinear elasticity, the described procedure can be used to provide an initial guess for Newton's method at the final loading step. Alternatively, $\f{u}_{lin}^N$ can also serve as a stand-alone displacement field defining the deformation $\pmb{\Phi}:\Omega_0\to\Omega$. This may be an interesting option since only a linear elasticity solver is required to implement it.  The linear diagonal incremental loading is also extensively used in ALE algorithms to deform the computational mesh for fluid domains in FSI problems \cite{stein2003mesh}. 

%% file: examples.tex
\section{Examples}
\subsection{2D single-patch domains}
First, we consider two two-dimensional, single-patch examples. We demonstrate the performance of the mesh deformation approach and show its dependence on the initial domain and on the value of Poisson's ratio used in the material law. As a rule, we use Newton's method with the Nonlinear Diagonal Incremental Loading (N-DIL). However, we also apply the Linear Diagonal Incremental Loading (L-DIL) as a stand-alone deformation method and compare the results. Finally, we compare the output of the mesh deformation approach with the results of the elliptic grid generation technique \cite{hinz2018elliptic} and the constrained optimization approach based on the area-orthogonality quality measure \cite{Gravesen2014}. When comparing different parametrizations, we mainly use the minimum of the Jacobian determinant
\begin{equation}\label{eq:minJ}
m(\f{G}) = \displaystyle\min_{\pmb{\xi}\in[0,1]^d}J(\f{G}(\pmb{\xi}))
\end{equation}
as the most neutral quality measure which does not favor any parametrization quality but its bijectivity. The higher the value of of $m(\f{G})$, the better. Secondary, we use the global ratio of the Jacobian determinant
\begin{equation}\label{eq:ratioJ}
R(\f{G}) = \frac{\displaystyle\max_{\pmb{\xi}\in[0,1]^d}J(\f{G}(\pmb{\xi}))}{\displaystyle\min_{\pmb{\xi}\in[0,1]^d}J(\f{G}(\pmb{\xi}))}
\end{equation}
as a measure of uniformity. The closer it is to 1, the better.
 
The mesh deformation approach is implemented using G+Smo \cite{jlmmz2014} - an open source C++ library providing necessary IGA routines. The area-orthogonality optimization is based on the nonlinear optimization library IPOPT. An in-house Newton-Krylov solver written in Python is used to implement the elliptic grid generation technique. 

\subsubsection{2D male rotor}
As the first example, we study the profile of a screw compressor's male rotor \cite{shamanskiy2018screw}. Its boundary is given as four cubic B-spline curves, and the domain is fairly simple so all considered parametrization techniques can be expected to perform well. We would like to notice, however, that the Coons patch does not produce a bijective parametrization when applied to this geometry.

\begin{figure}[H]
	\centering
	\includegraphics[width=0.49\textwidth]{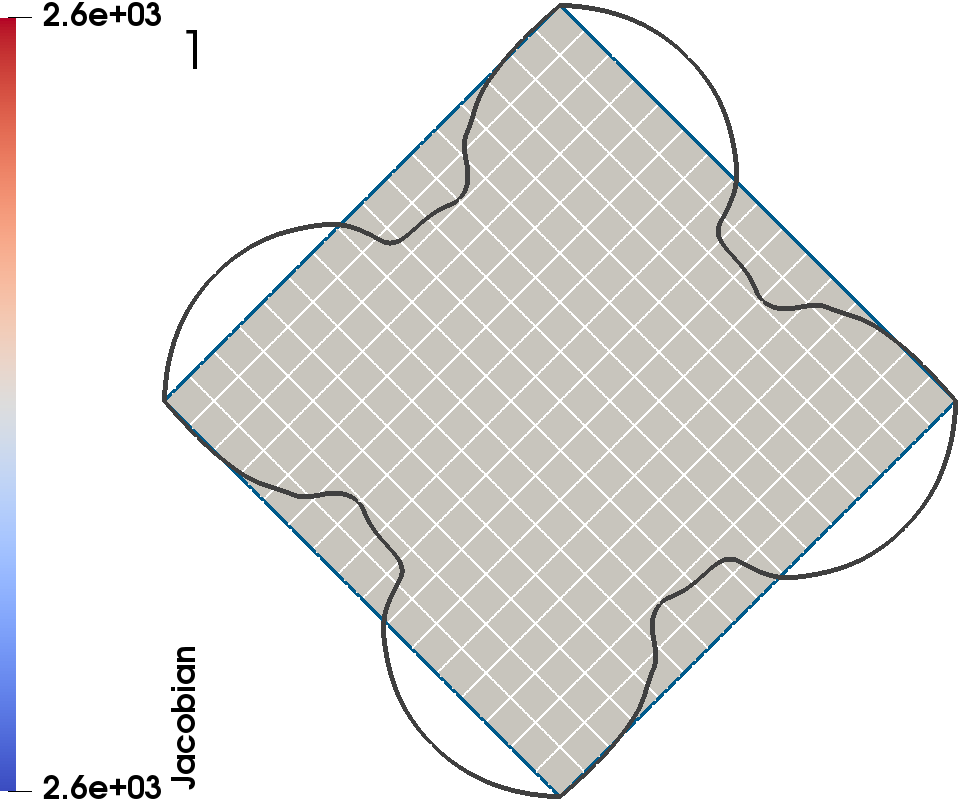}
	\includegraphics[width=0.49\textwidth]{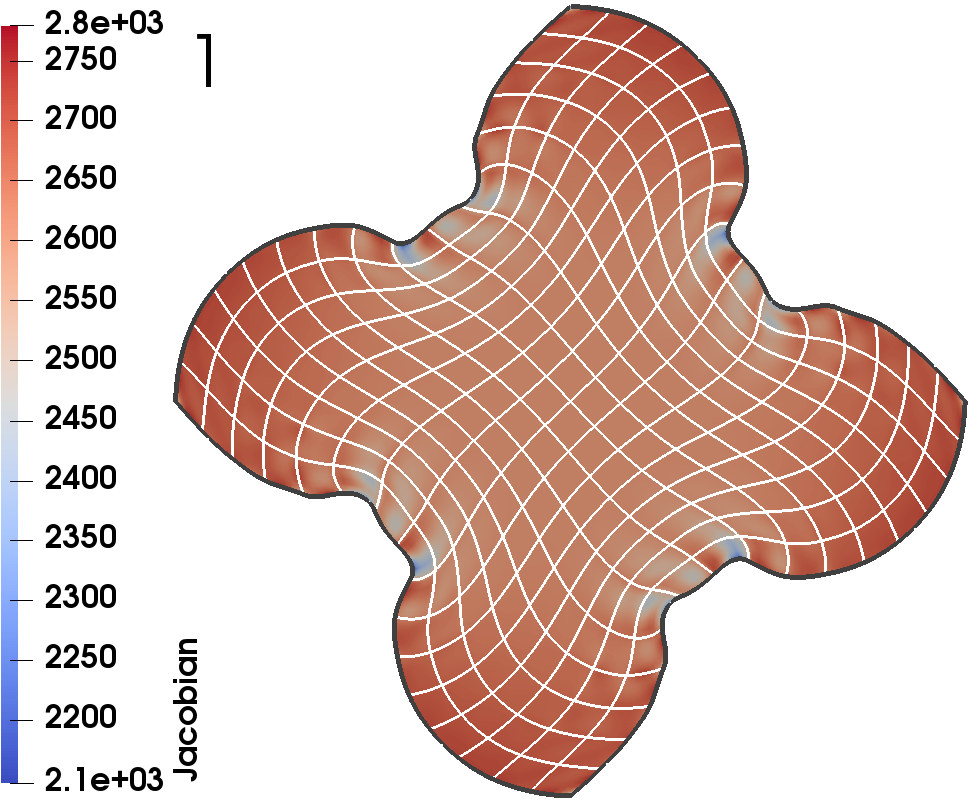}
	\includegraphics[width=0.49\textwidth]{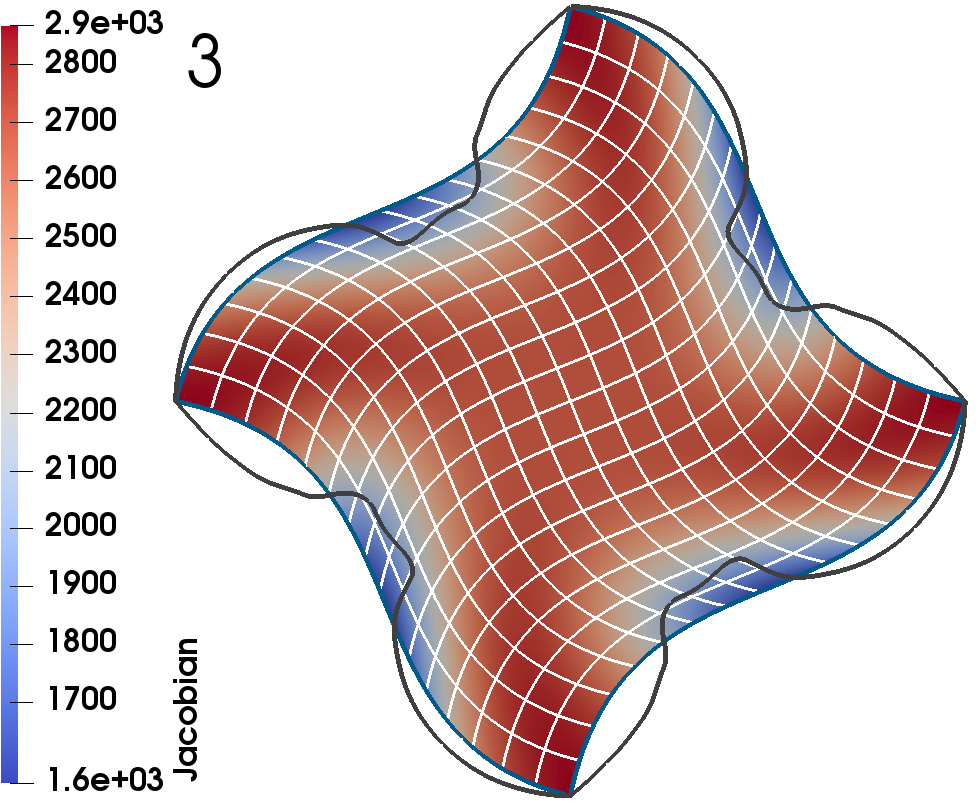}
	\includegraphics[width=0.49\textwidth]{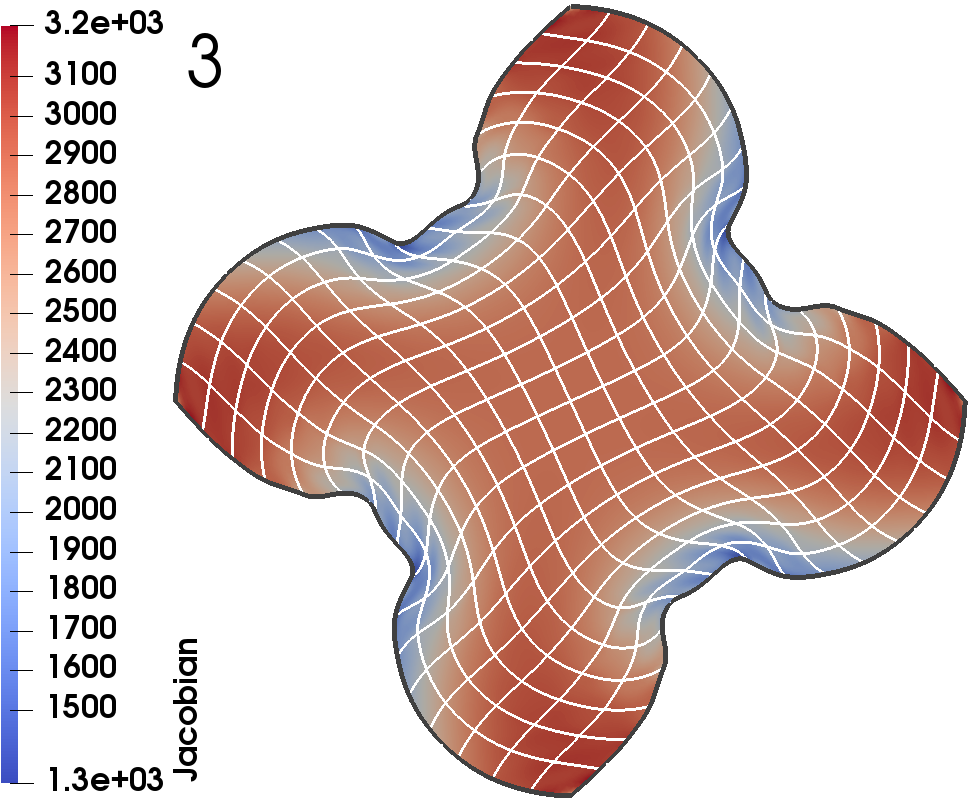}
	\caption{Male rotor example. Initial domains ($p = 1$ and $p=3$) and results of deformation by Newton's method with N-DIL and Poissson's ratio of 0.49.}\label{fig:male}
\end{figure}

Figure~\ref{fig:male} depicts the results of computing the deformation using Newton's method with N-DIL. Two initial domains were generated by applying $L^2$-simplification to each part of the boundary with the coarsest B-Spline bases of degree $p=1$ and $p=3$. A value of 0.49 was used for Poisson's ratio. Such a high value required us to use at least $N=5$ loading steps for the initial domain with $p=1$ and $N=3$ for the initial domain with $p=3$. Both initial domains led to high-quality bijective parametrizations. With respect to the quality measures $m(\f{G})$ and $R(\f{G})$, the $p=1$ parametrization is better. It inherited its uniform structure from the initial domain due to the high value of Poisson's ratio. We use it a as baseline for the following comparison.

\begin{figure}[H]
	\centering
	\includegraphics[height=6cm]{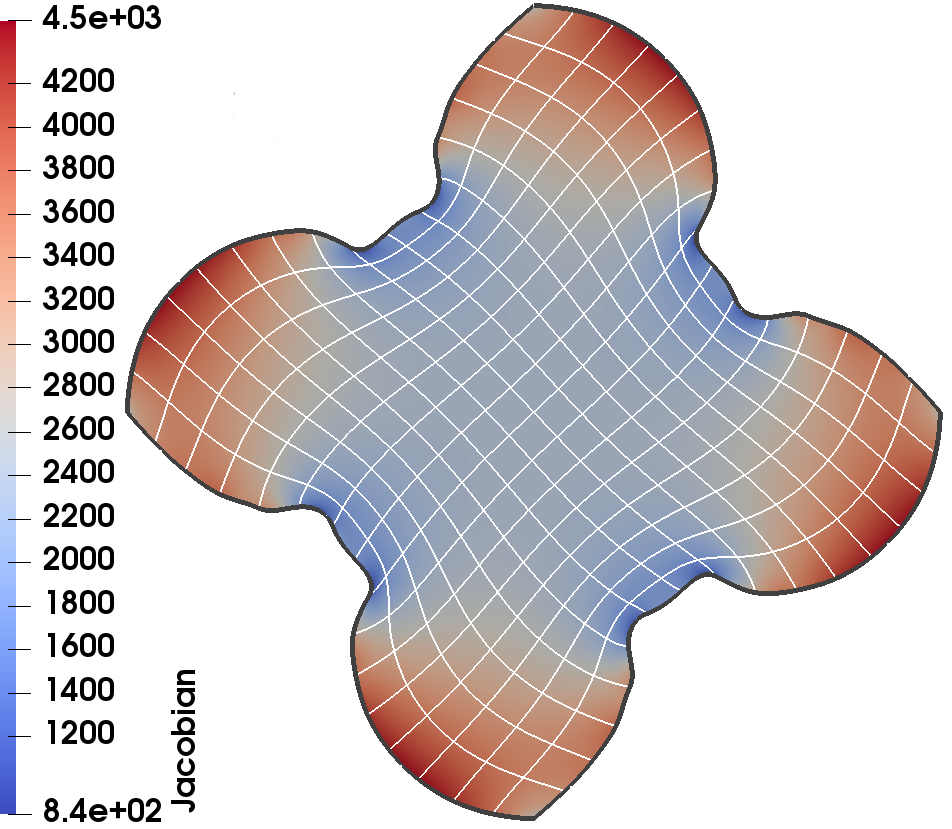}
	\includegraphics[height=6cm]{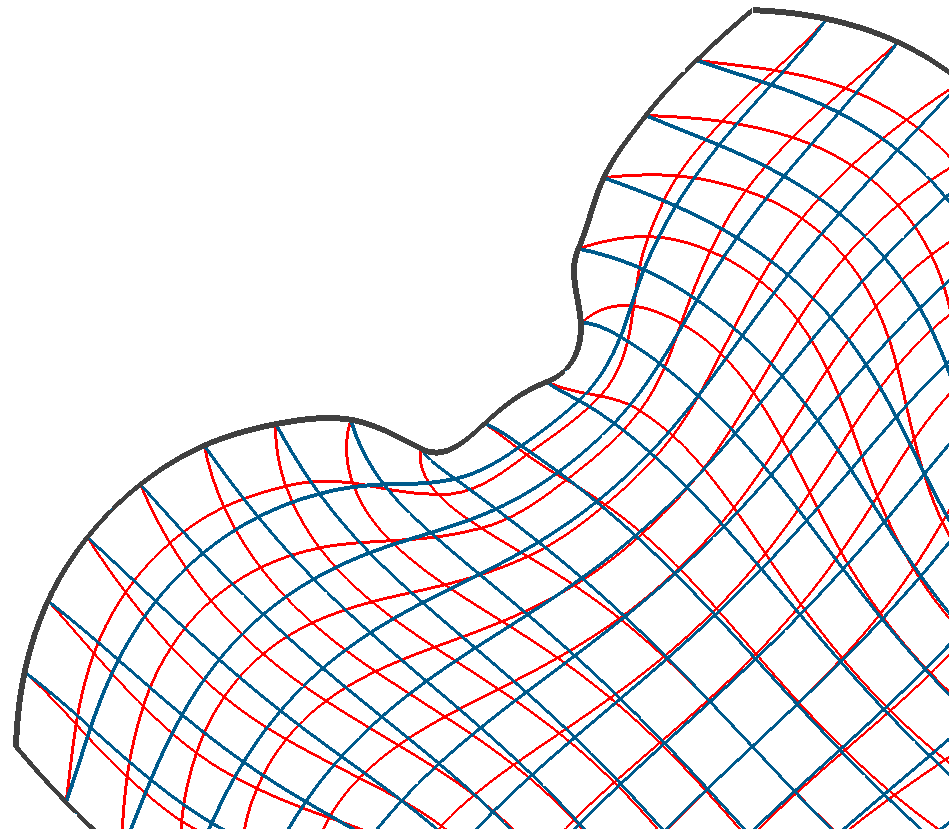}
	\caption{Male rotor example. Initial domain $p=1$ deformed by Newton's method with N-DIL and Poisson's ratio of 0 (left). Comparison of the blue corresponding mesh with the red baseline mesh (right).}\label{fig:malePoiss}
\end{figure}

The result of the mesh deformation approach depends heavily on the choice of Poisson's ratio. To illustrate it, we applied Newton's method with N-DIL to deform the $p=1$ initial domain with Poisson's ratio of 0, see Fig.~\ref{fig:malePoiss}. The resulting parametrization is still bijective but is worse than the baseline with respect to $m(\f{G})$ and $R(\f{G})$. However, a lower value of Poisson's ratio made it easier to compute a bijective displacement field; only $N=1$ loading step was necessary.

\begin{figure}[H]
	\centering
	\includegraphics[height=6cm]{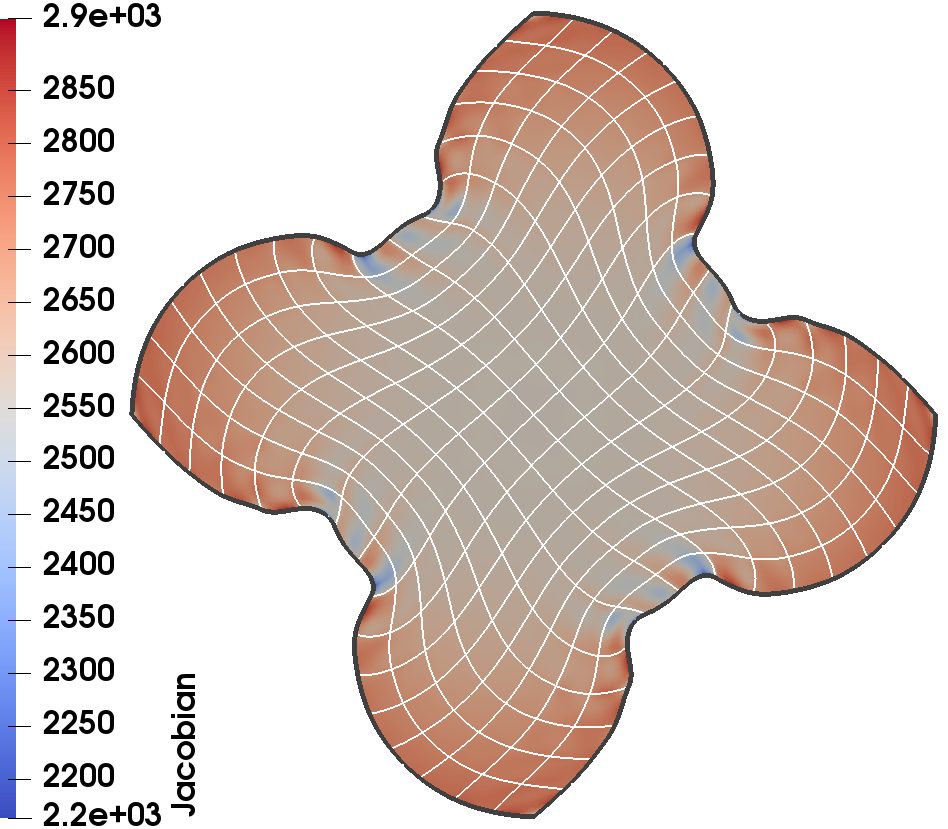}
	\includegraphics[height=6cm]{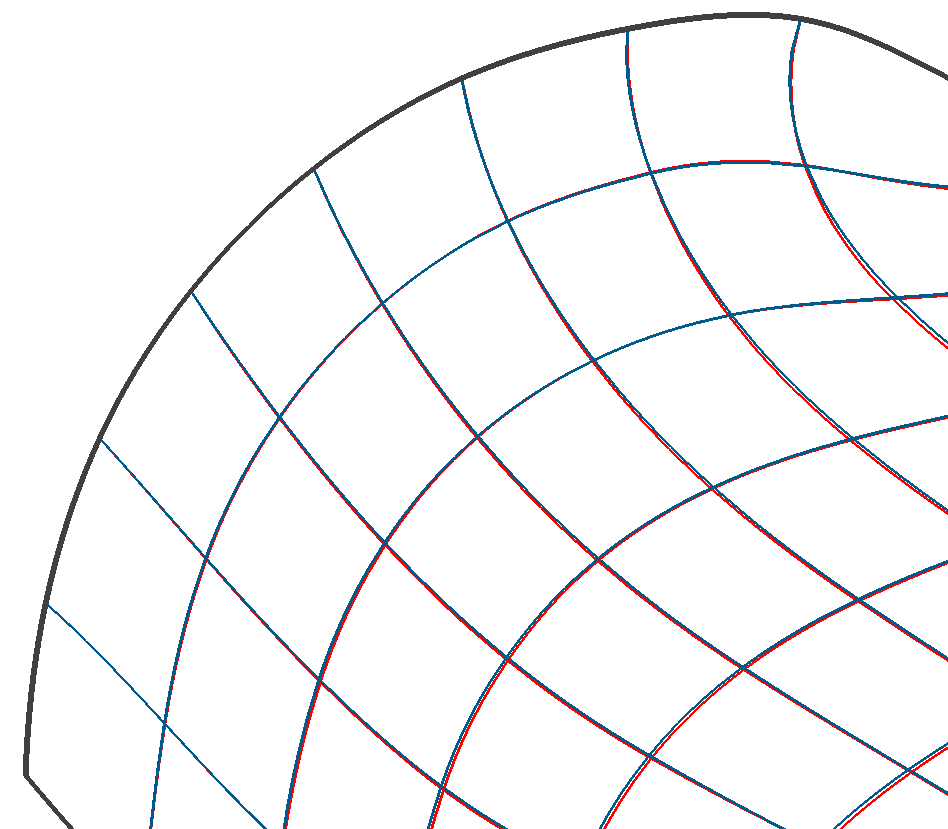}
	\caption{Male rotor example. Initial domain $p=1$ deformed by the L-DIL method with $N=10$ loading steps and Poisson's ratio of 0.49 (left). Comparison of the blue corresponding mesh with the red baseline mesh (right).}\label{fig:maleLin}
\end{figure}

Additionally, we show the performance of the L-DIL method as a stand-alone deformation technique. Figure~\ref{fig:maleLin} presents the results of applying it to deform the $p=1$ initial domain with Poisson's ratio of 0.49. After $N=10$ loading steps, the resulting parametrization is virtually indiscernible from the baseline. In order to achieve bijectvity, at least $N=3$ loading steps had to be used.

\begin{figure}[H]
	\centering
	\includegraphics[height=6cm]{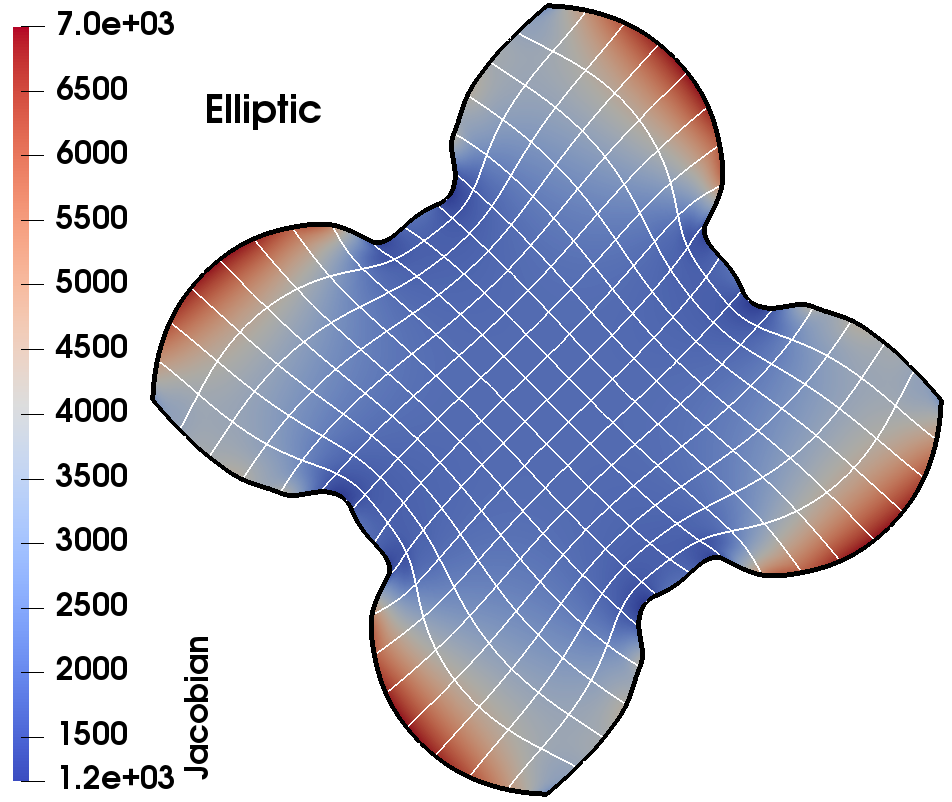}
	\includegraphics[height=6cm]{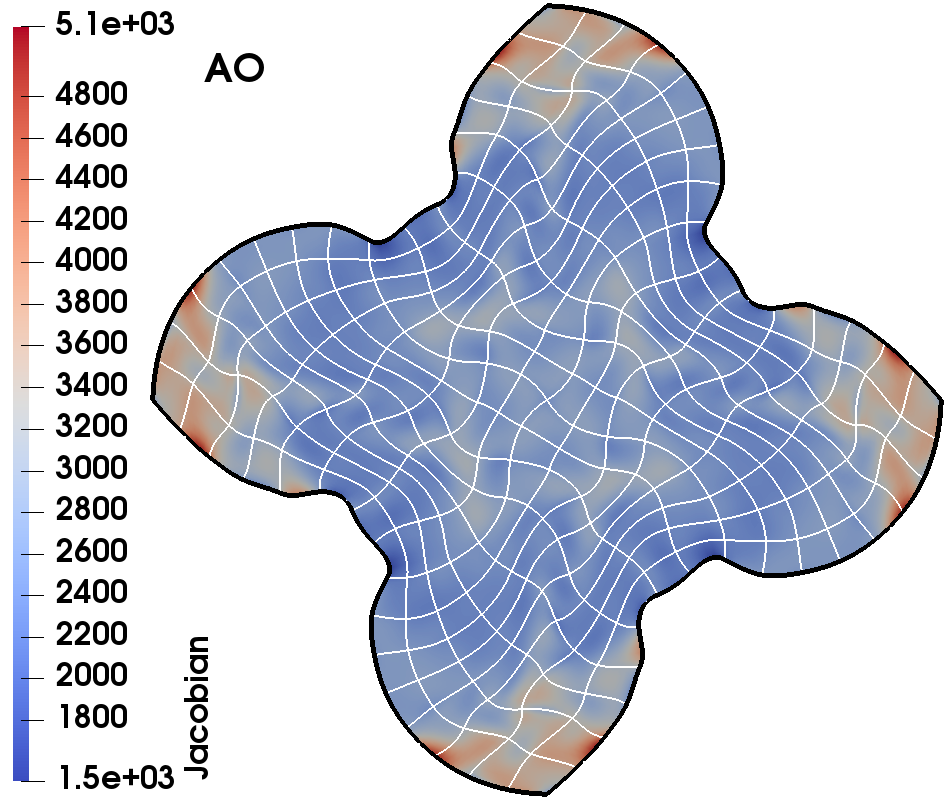}
	\caption{Male rotor example. Parametrizations by elliptic grid generation (left) and area-orthogonality optimization (right).}\label{fig:maleOther}
\end{figure}

Finally, we applied elliptic grid generation and area-orthogonality optimization to enrich the comparison, see Fig.~\ref{fig:maleOther}. Both techniques produce high-quality bijective parametrizations; however, the baseline is better with respect to $m(\f{G})$ and $R(\f{G})$.

\subsubsection{2D puzzle piece}
Next, we consider a puzzle piece example. Its boundary possesses distinct protruding and concave regions which make it difficult to construct a bijective tensor-product parametrization.

\begin{figure}[H]
	\centering
	\includegraphics[width=0.49\textwidth]{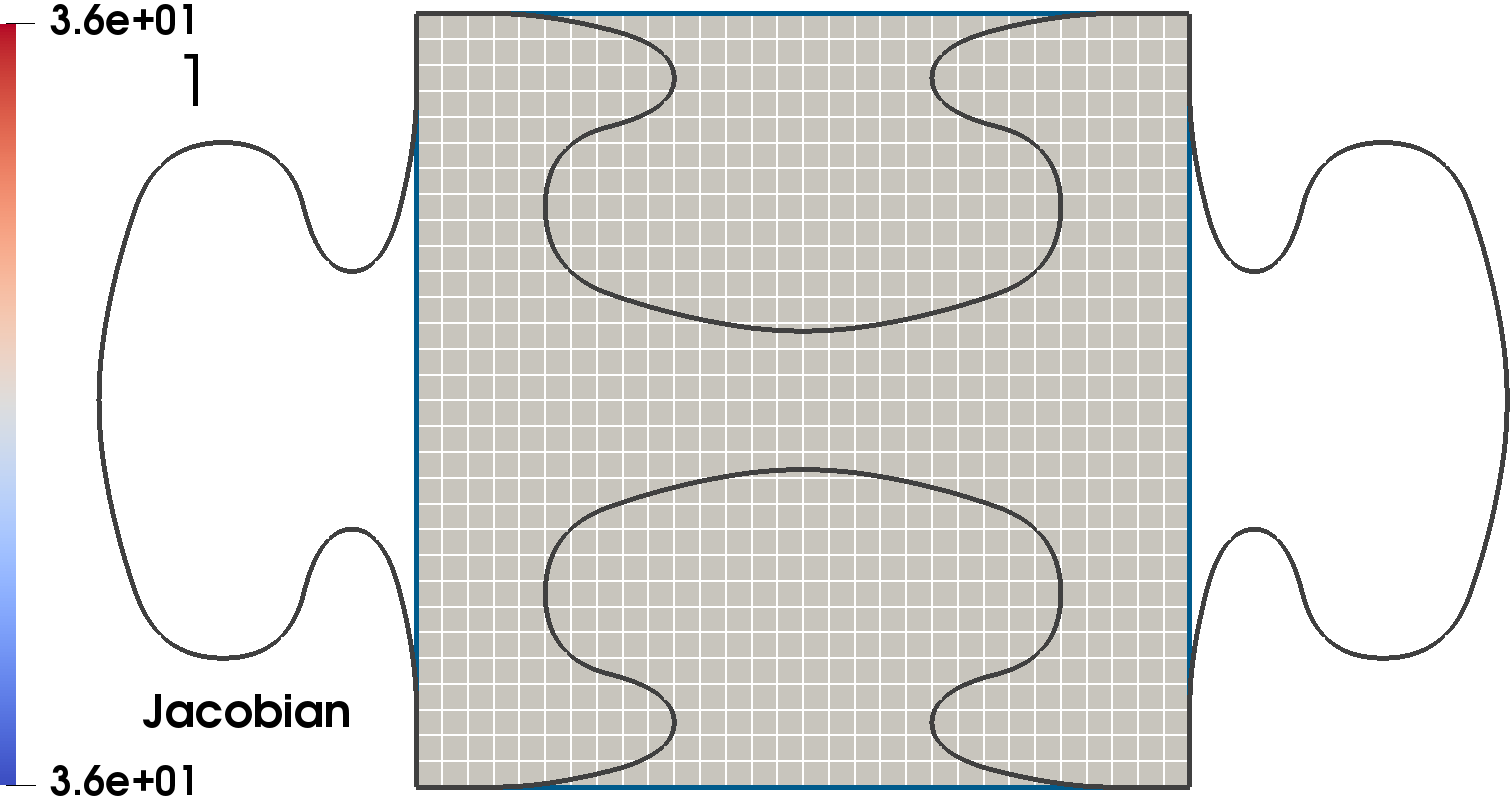}
	\includegraphics[width=0.49\textwidth]{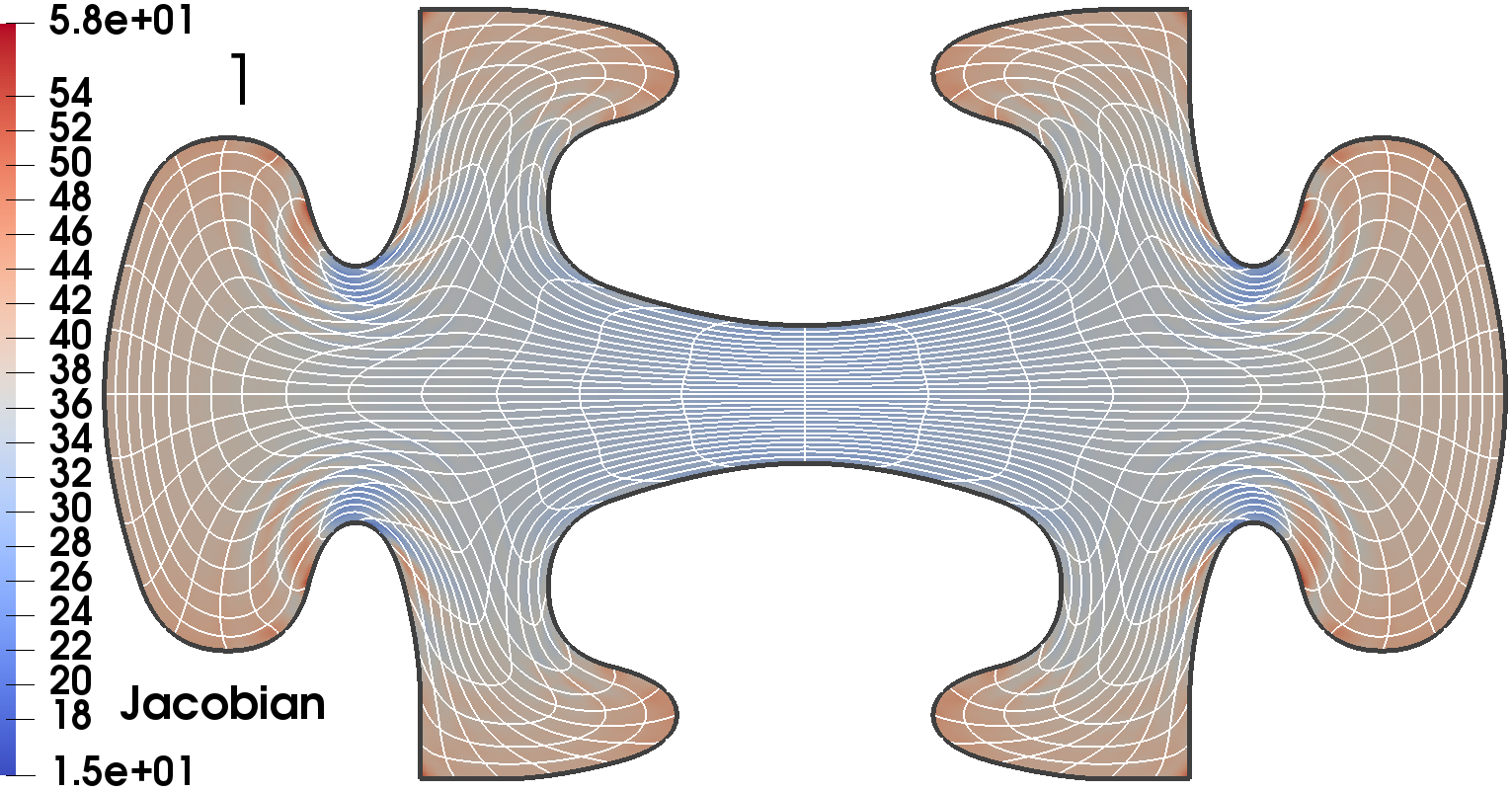}\\
	\includegraphics[width=0.49\textwidth]{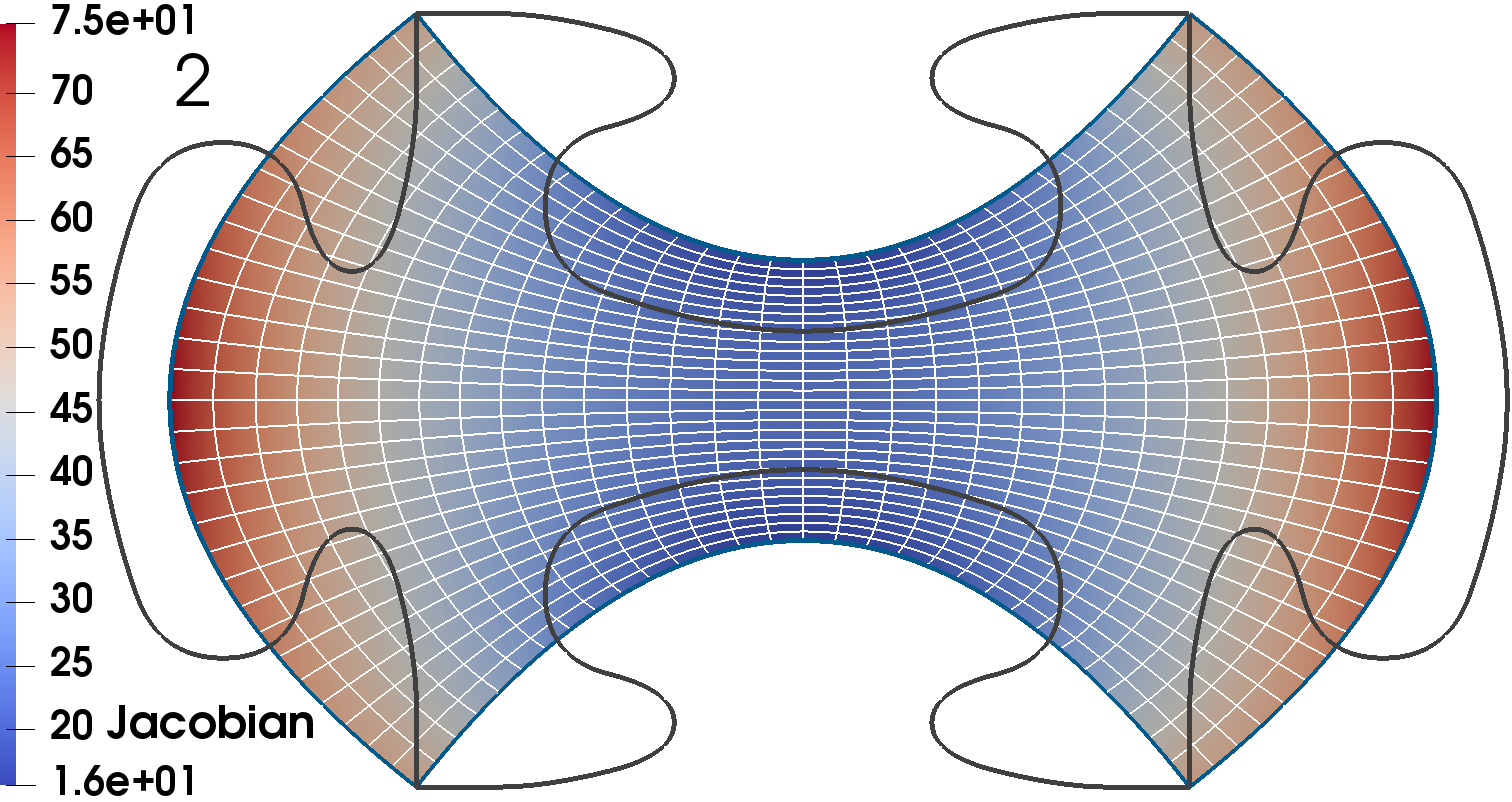}
	\includegraphics[width=0.49\textwidth]{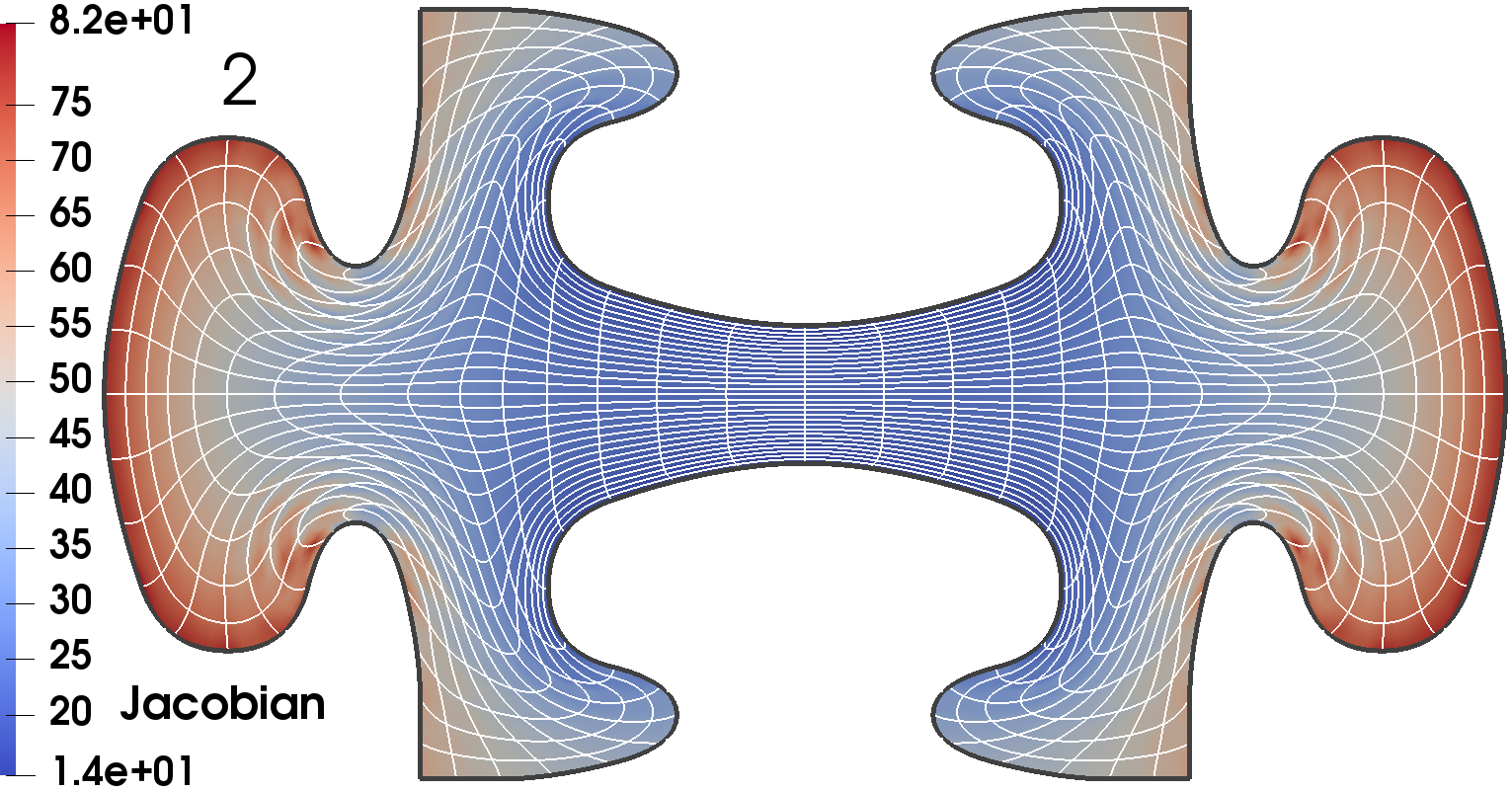}\\
	\caption{Puzzle piece example. Initial domains ($p = 1$ and $p=2$) and results of deformation by Newton's method with N-DIL and Poissson's ratio of 0.49.}\label{fig:puzzle}
\end{figure}

We generated two different initial domains by applying the $L^2$-simplification to each part of the boundary with the coarsest B-spline bases of degree $p=1$ and $p=2$. Newton's method with N-DIL was applied with Poisson's ratio equal to 0.49. Such a high value, together with the complexity of the domain, made it necessary to use $N=13$ loading steps for the $p=1$ initial domain and $N=8$ for the $p=2$ initial domain. Figure~\ref{fig:puzzle} depicts the results of the deformation. Judging by $m(\f{G})$ and $R(\f{G})$, the $p=1$ initial domain results in a better parametrization. However, the middle neck-like region of the domain underwent a large deformation which resulted in the isoparametric lines being pushed away to the sides. On the other hand, the $p=2$ initial domain is geometrically much closer to the target domain so it had to be deformed less. This results in a visually more natural parametrization which we use as a baseline for the following comparison.

\begin{figure}[H]
	\centering
	\includegraphics[width=0.49\textwidth]{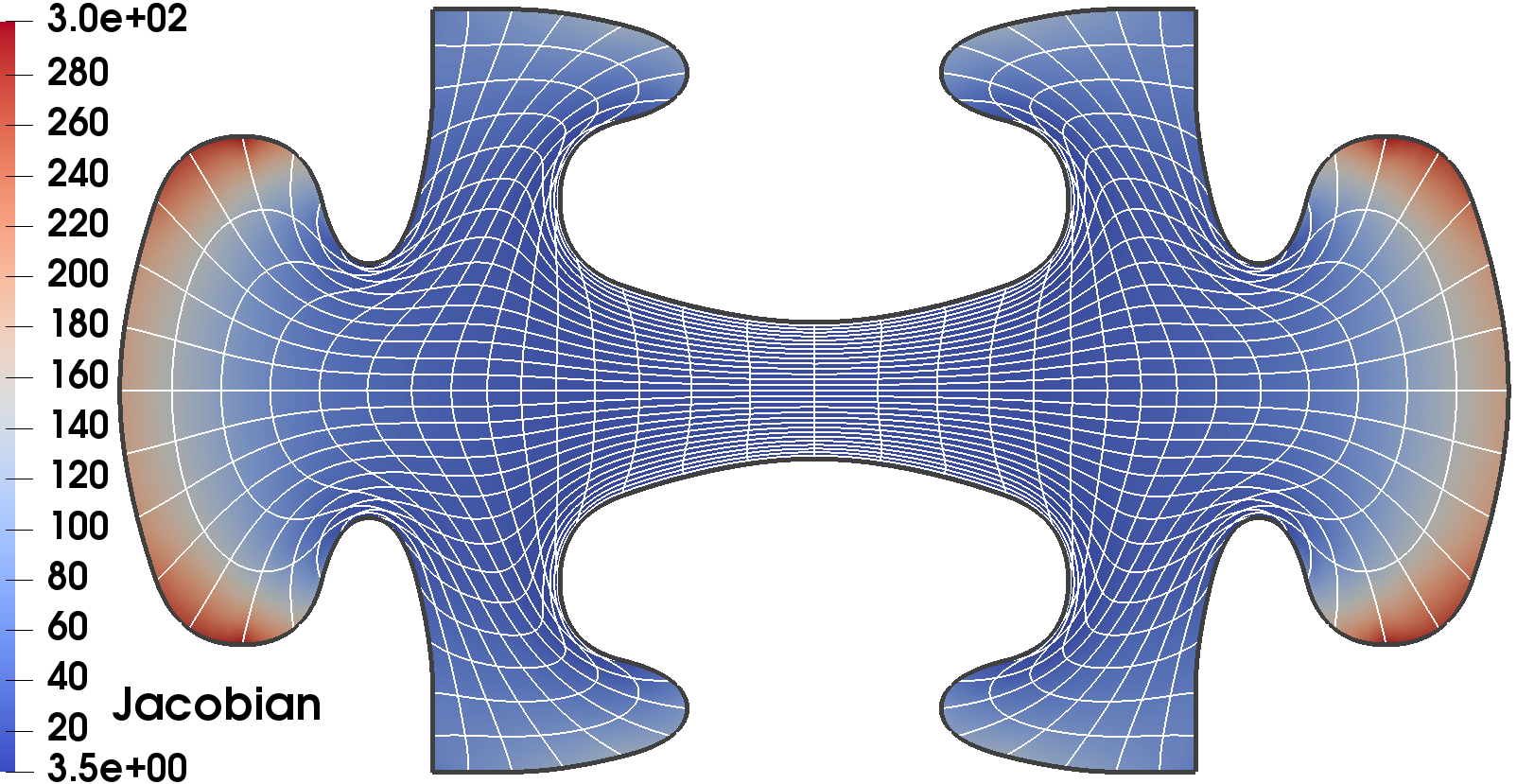}
	\includegraphics[width=0.49\textwidth]{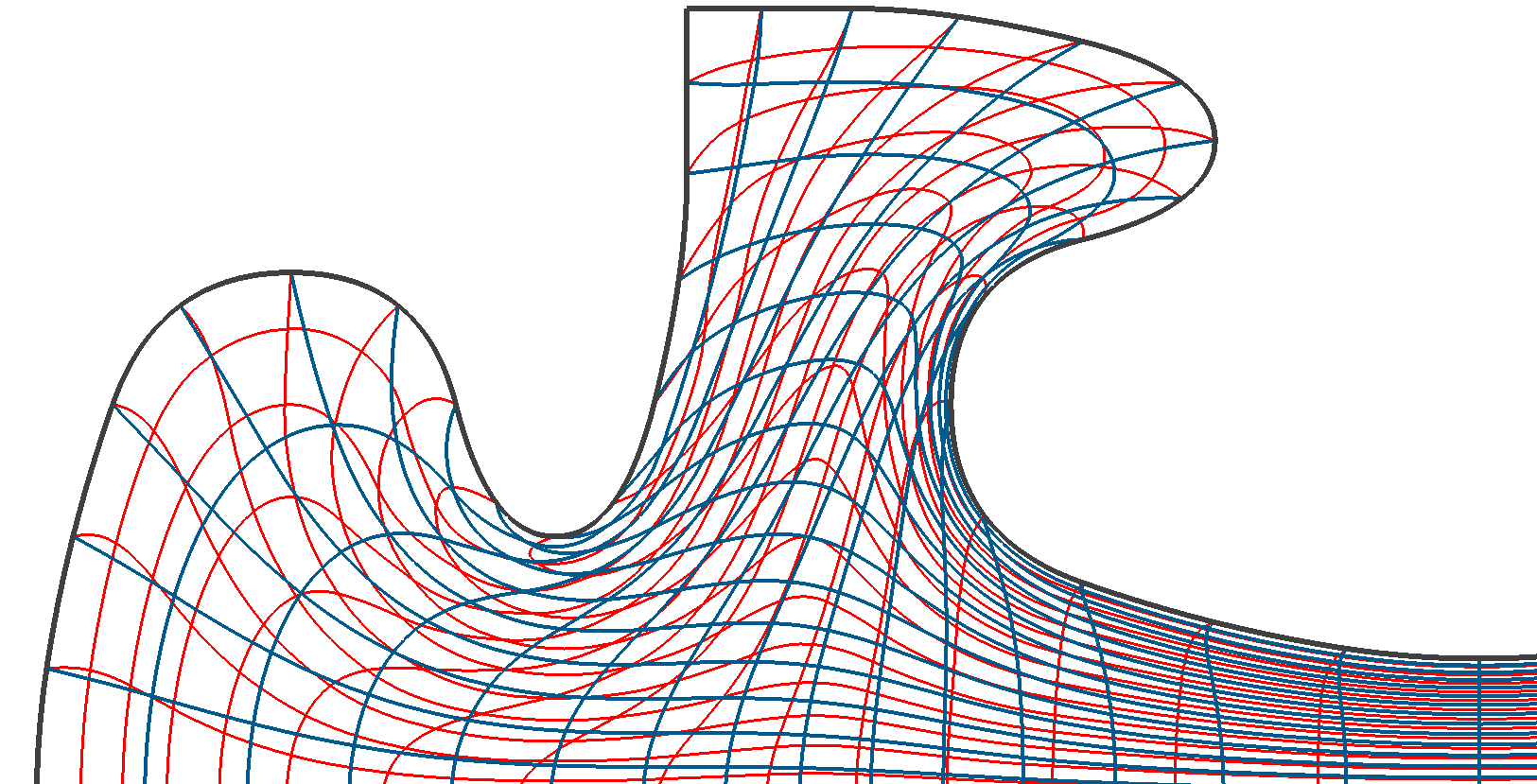}\\
	\caption{Puzzle piece example. Initial domain $p=2$ deformed by Newton's method with N-DIL and Poisson's ratio of 0 (left). Comparison of the blue corresponding mesh with the red baseline bash (right).}\label{fig:puzzlePoiss}
\end{figure}

Due to the complexity of the domain, it is crucial to use a high value of Poisson's ratio to preserve bijectivity. Figure~\ref{fig:puzzlePoiss} demonstrates the results of deforming the $p=2$ initial domain by Newton's method with N-DIL and Poisson's ratio equal to 0. The isoparametric lines come together densely next to the concave parts of the boundary, and $m(\f{G})$ drops almost by one order of magnitude in comparison to the baseline making the parametrization almost not bijective.  

\begin{figure}[H]
	\centering
	\includegraphics[width=0.49\textwidth]{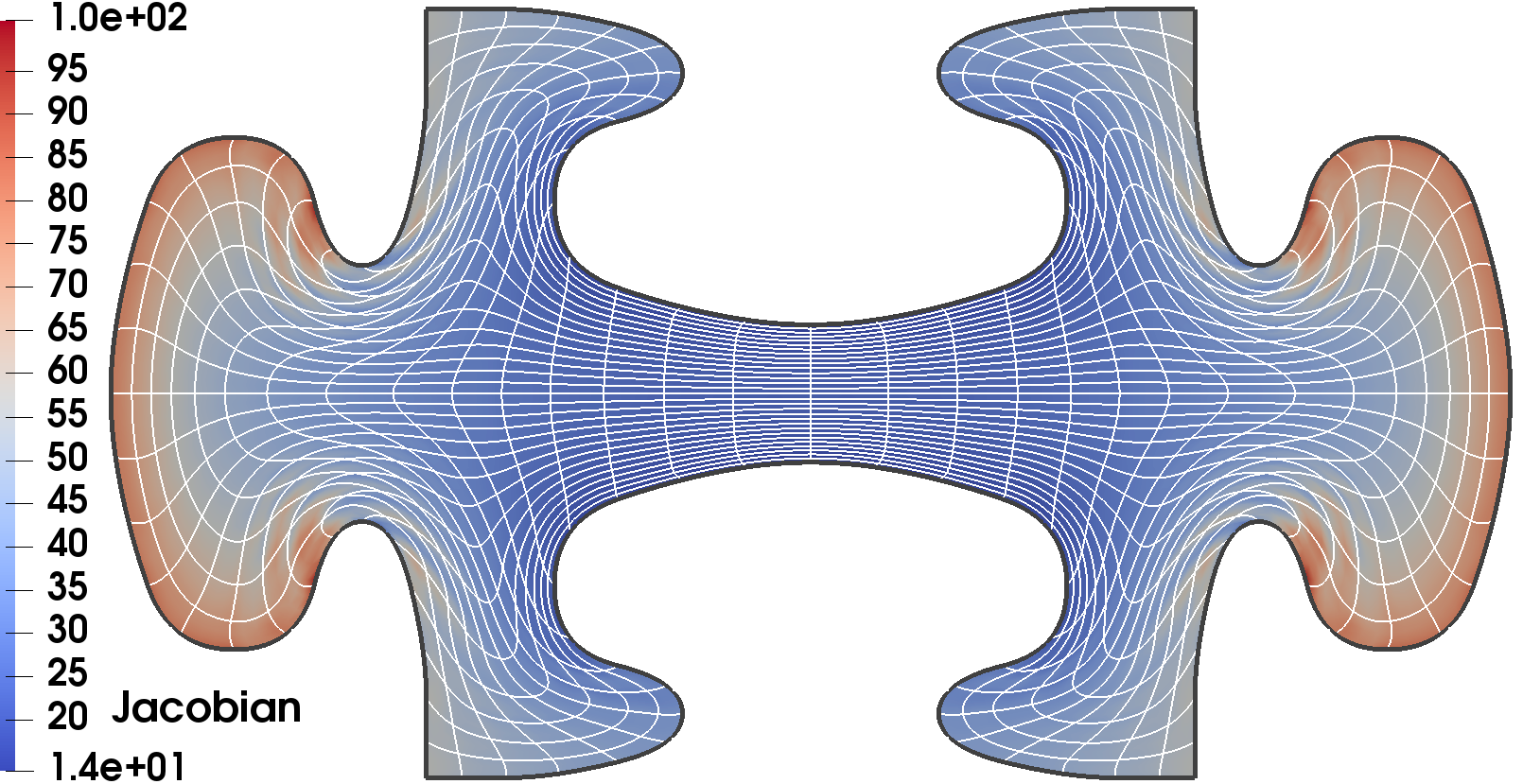}
	\includegraphics[width=0.49\textwidth]{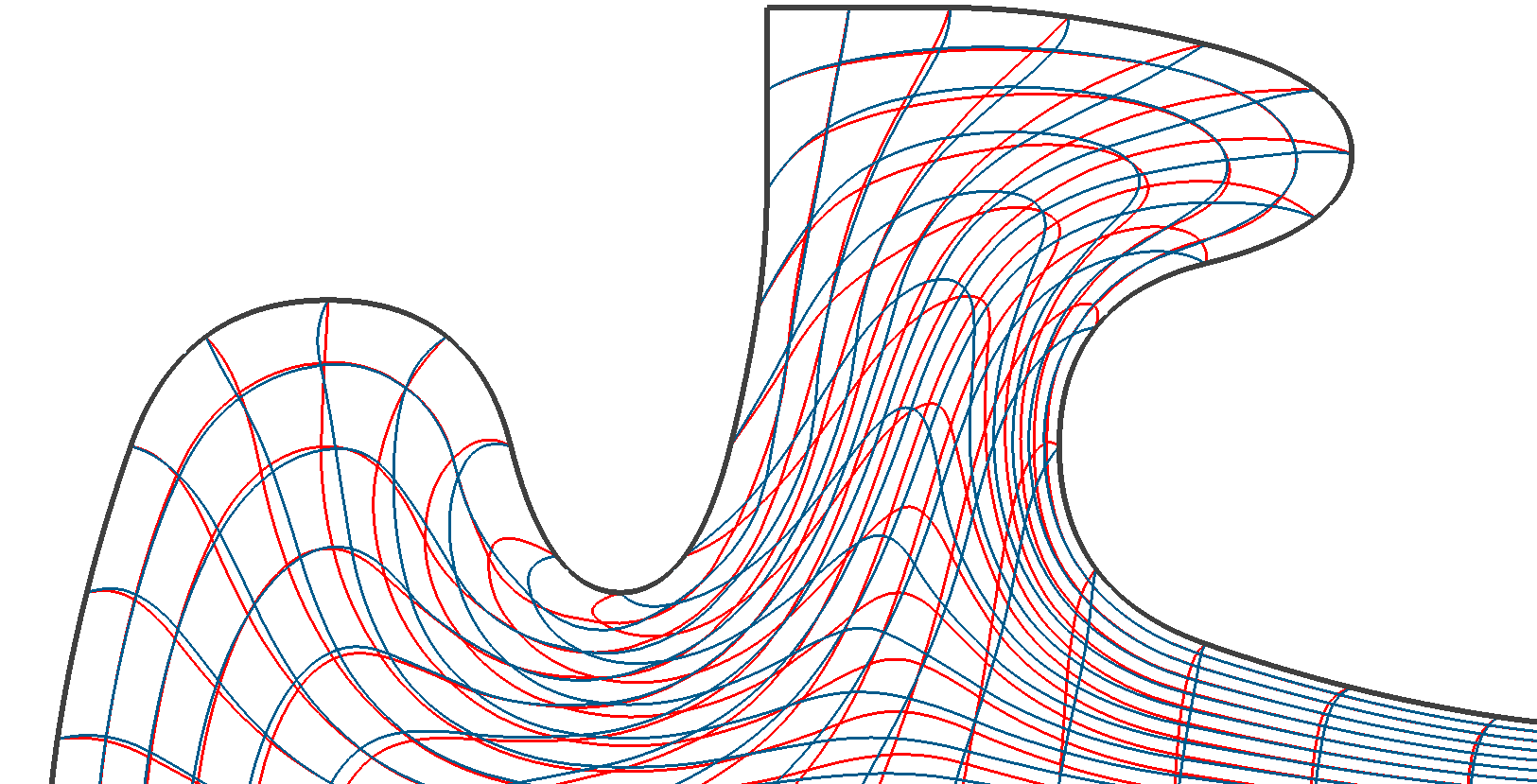}\\
	\caption{Puzzle piece example. Initial domain $p=2$ deformed by the L-DIL method with $N=10$ loading steps and Poisson's ratio of 0.49 (left). Comparison of the blue corresponding mesh with the red baseline mesh (right).}\label{fig:puzzleLin}
\end{figure}

Additionally, we demonstrate the performance of the L-DIL method as a stand-alone parametrization technique in Figure~\ref{fig:puzzleLin}. The $p=2$ initial domain was deformed with Poisson's ratio of 0.49 and $N=15$ loading steps. Unlike in the male rotor example, the resulting parametrization is quite different from the baseline. Still, it is bijective and has the same value of $m(\f{G})$; however, the baseline is more uniform. At least $N=8$ loading step are required to achieve bijectivity.

Finally, we applied elliptic grid generation and area-orthogonality optimization to the puzzle piece example. The former provides a barely bijective, highly non-uniform parametrization. The latter provides a high-quality parametrization, only slightly worse with respect to $m(\f{G})$ than the baseline.

\begin{figure}[H]
	\centering
	\includegraphics[width=0.49\textwidth]{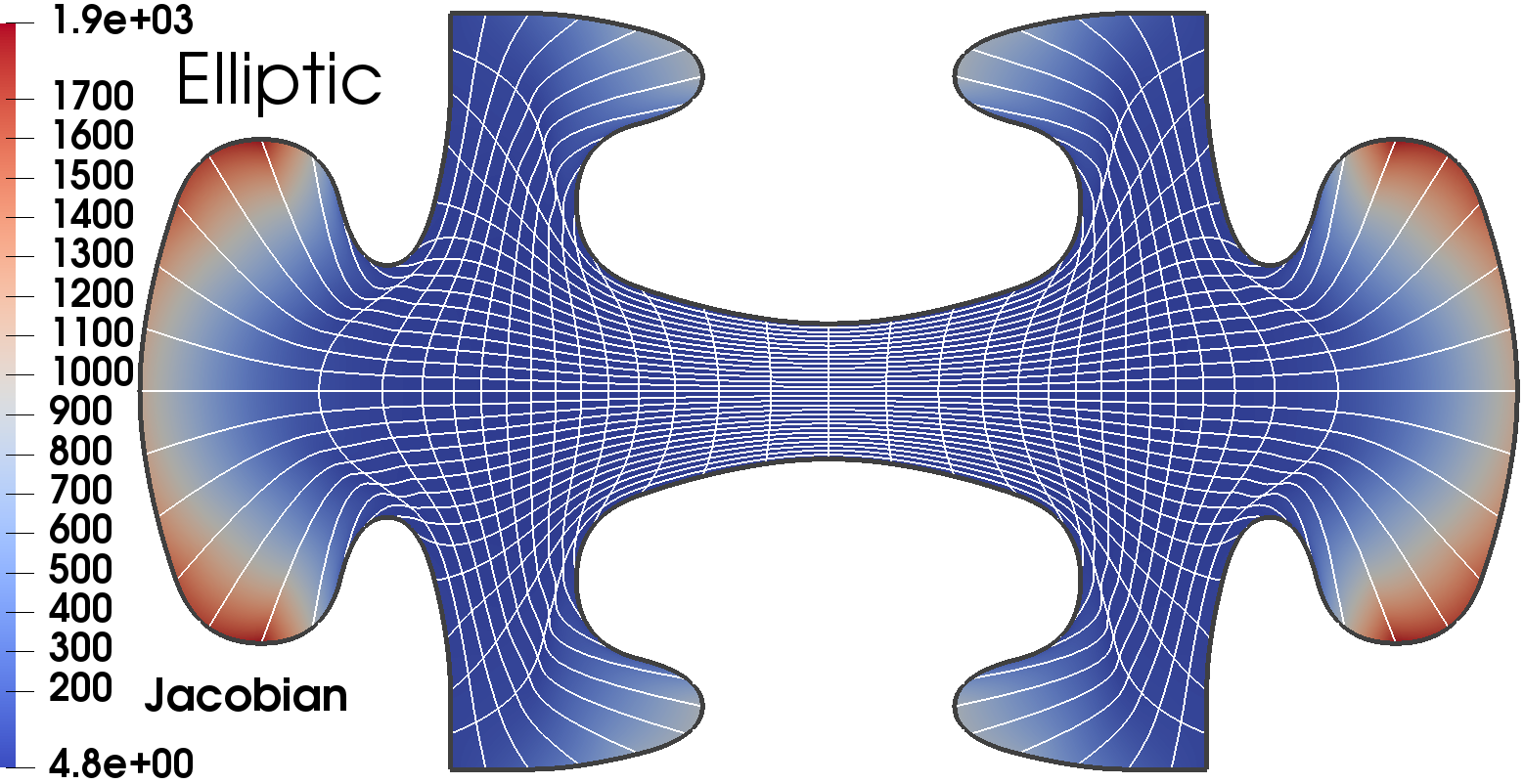}
	\includegraphics[width=0.49\textwidth]{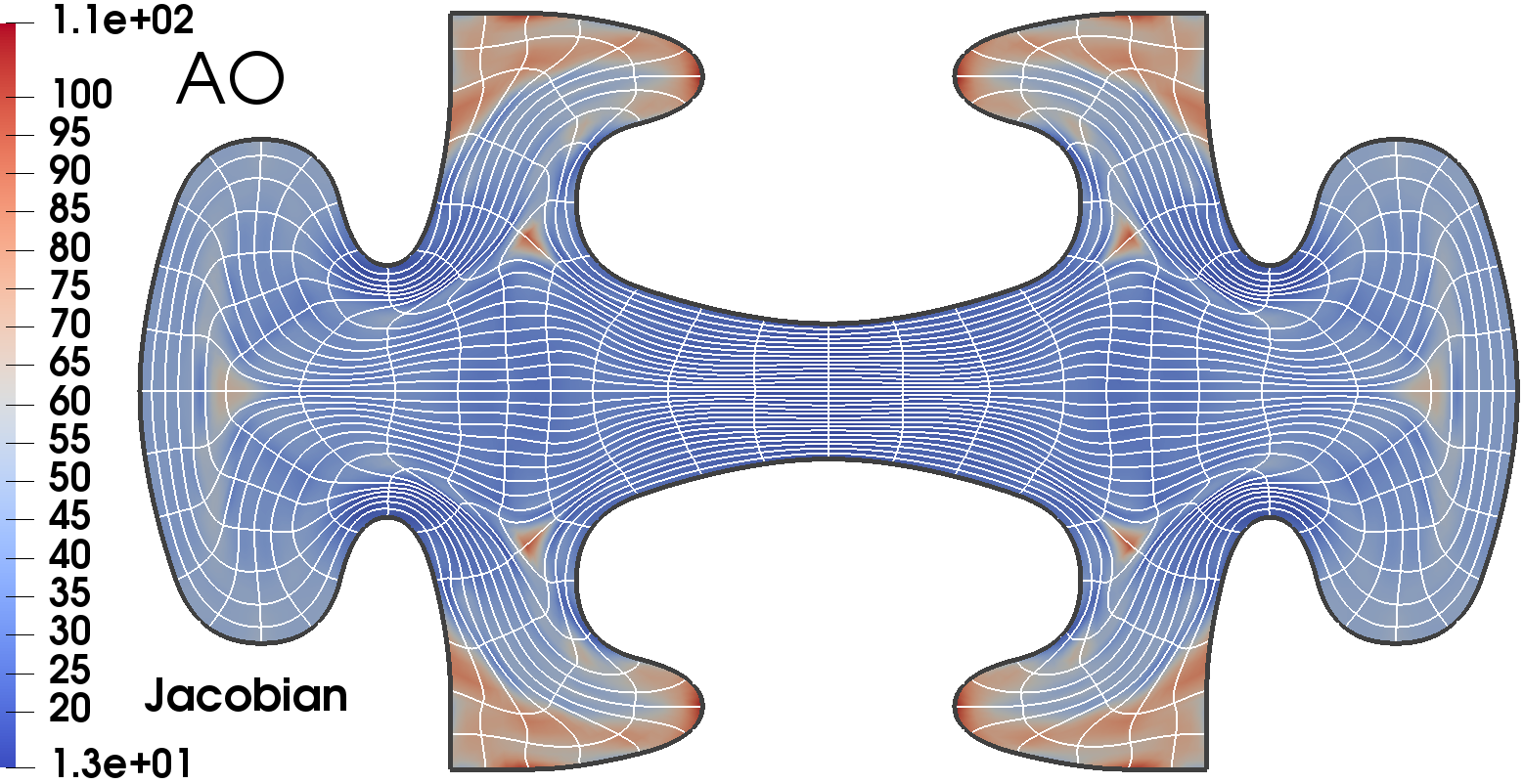}\\
	\caption{Puzzle piece example. Parametrizations by elliptic grid generation (left) and area-orthogonality optimization (right).}\label{fig:puzzleOther}
\end{figure}

\subsubsection*{Remark on numerical effort}
Here we briefly describe our experience with respect to the numerical cost of the applied parametrization approaches. Unfortunately, since they are implemented in different programming languages, a fair comparison of CPU time necessary for every method to produce a bijective parametrization is not possible. We can, however, get an impression of their numerical cost by looking at the number of iterations taken by each method. In our experience, elliptic grid generation is the fastest method which takes only 3-6 iterations to converge. When applying the mesh deformation approach,  3-10 loading steps are required to produce a bijective initial guess by the diagonal incremental loading. The result can be used as a final parametrization, or additional 4-7 iterations of Newton's method are necessary to acquire a solution to the system (\ref{eq:nonlin}-\ref{eq:DBC}). Together, this results in 7-17 iterations. Finally, the optimization technique takes 40-70 iterations which makes it the most computationally expensive.

\subsubsection*{Convergence of diagonal incremental loading}
\begin{figure}[H]
	\centering
	\includegraphics[height = 6cm]{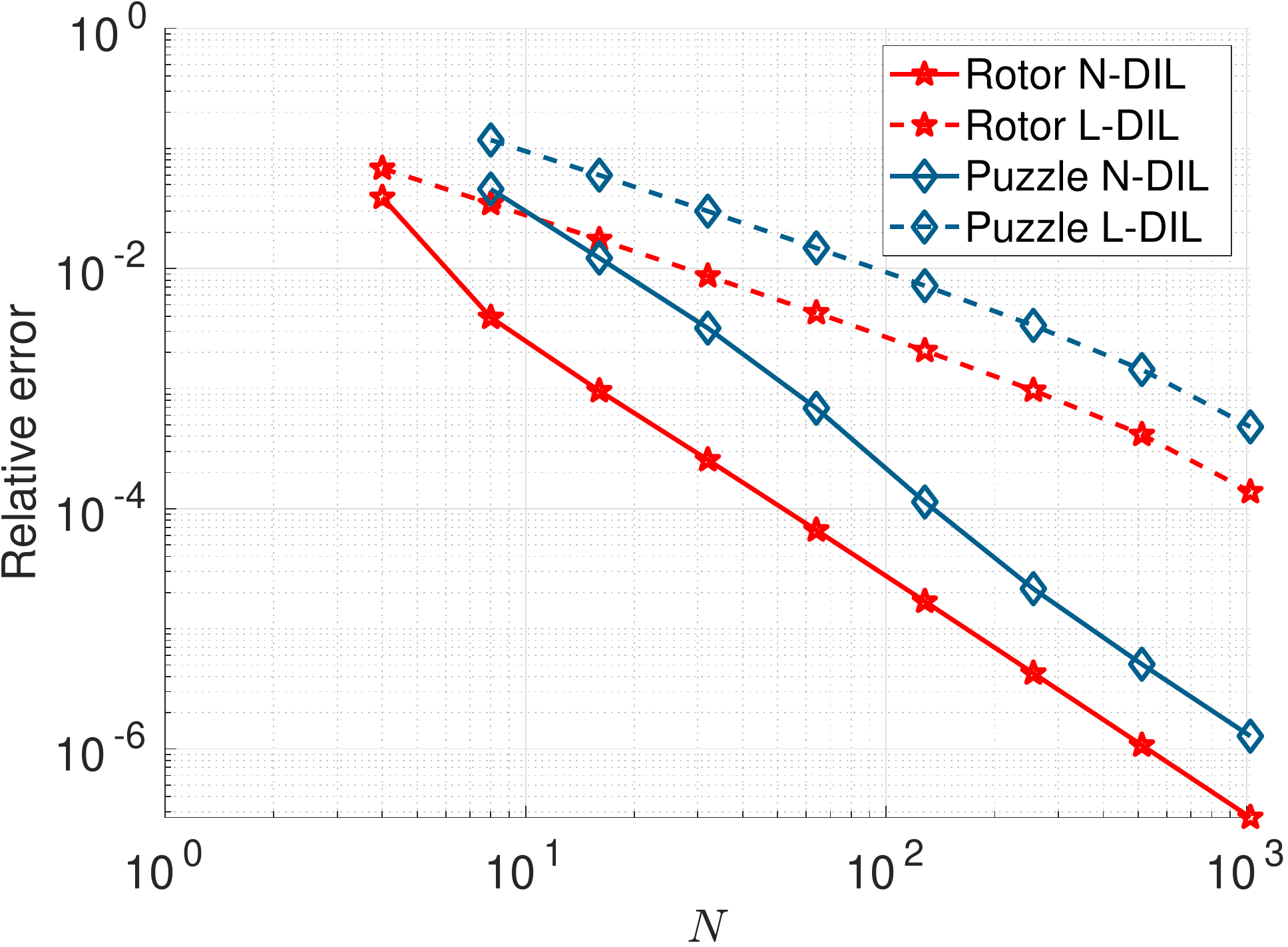}
	\caption{Convergence of the nonlinear and linear diagonal incremental loading algorithms for the male rotor and the puzzle piece examples. }\label{fig:conv}
\end{figure}

We conclude the analysis of different aspects of the mesh deformation approach by studying the convergence of the nonlinear and linear diagonal incremental loading approaches. As we mention in Section 5, the result of N-DIL $\f{u}^N_{inc}$ converges quadratically to the solution of the system (\ref{eq:nonlin}-\ref{eq:DBC}) $\f{u}$ as the number of loading steps $N$ grows. At the same time, the result of L-DIL $\f{u}_{lin}^N$ converges linearly to a different displacement $\f{u}_{lin}$ which is quiet close to $\f{u}$ and, surprisingly, bijective. Figure~\ref{fig:conv} presents a convergence plot where the relative errors
\begin{equation}
err_{N-DIL} = \frac{||\f{u}-\f{u}^N_{inc}||_{L^2}}{||\f{u}||_{L^2}} \; \text{ and }\; err_{L-DIL} = \frac{||\f{u}_{lin}-\f{u}^N_{lin}||_{L^2}}{||\f{u}_{lin}||_{L^2}} 
\end{equation}
are plotted against $N$ for the both examples.

\subsection{2D multi-patch female rotor}
\begin{figure}[H]
	\centering
	\includegraphics[height=6.1cm]{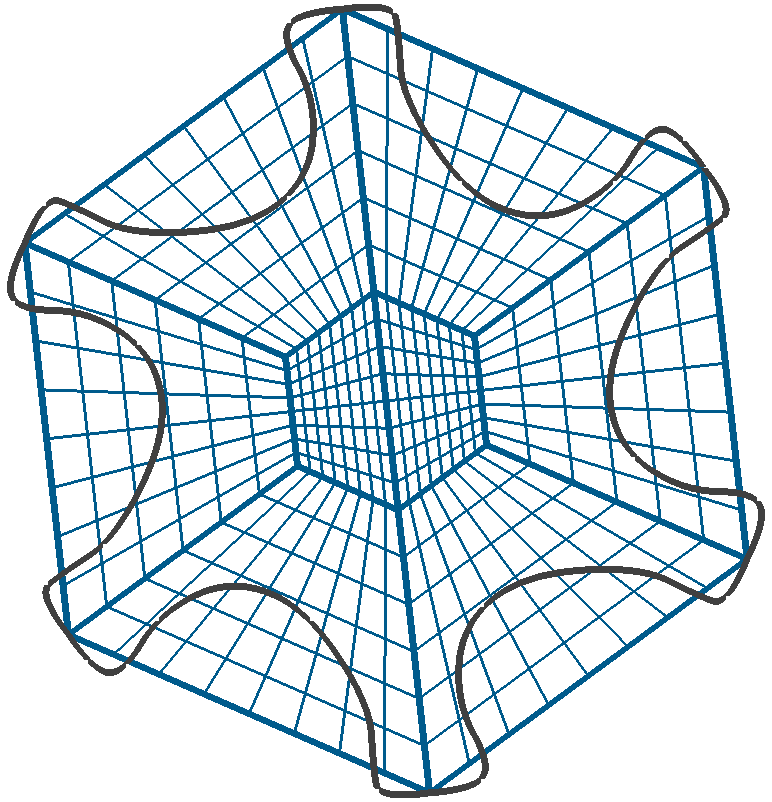}
	\includegraphics[height=6.1cm]{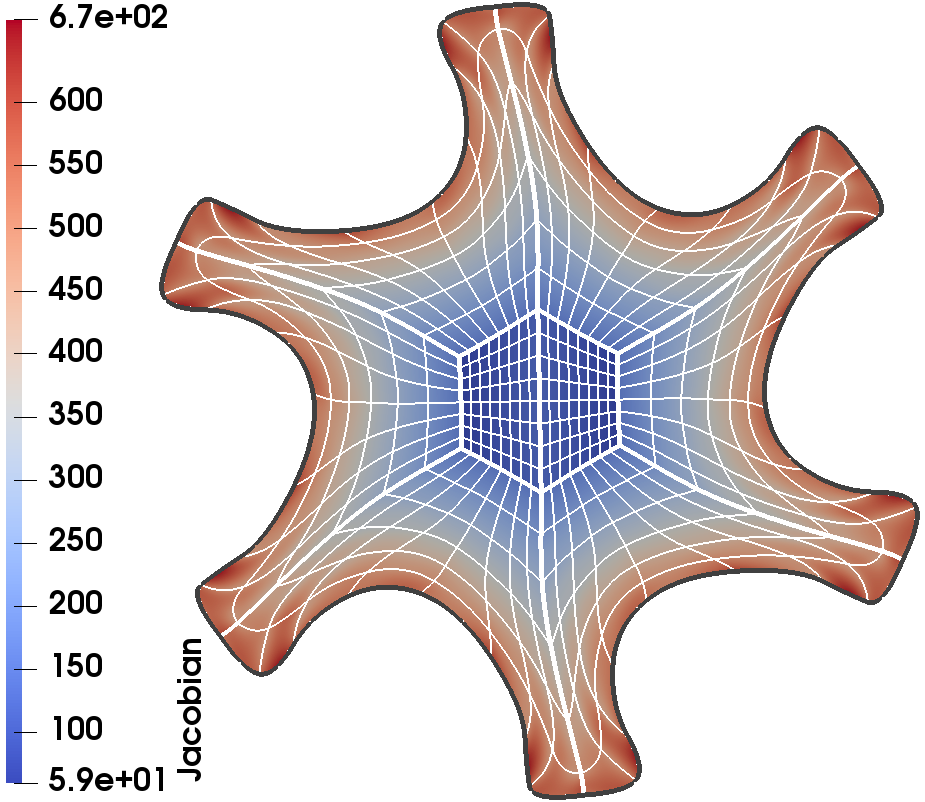}
	\caption{Mesh deformation approach for the female rotor. Initial domain (left) and resulting parametrization (right).}\label{fig:female}
\end{figure}

\begin{figure}[H]
	\centering
	\includegraphics[height=6.1cm]{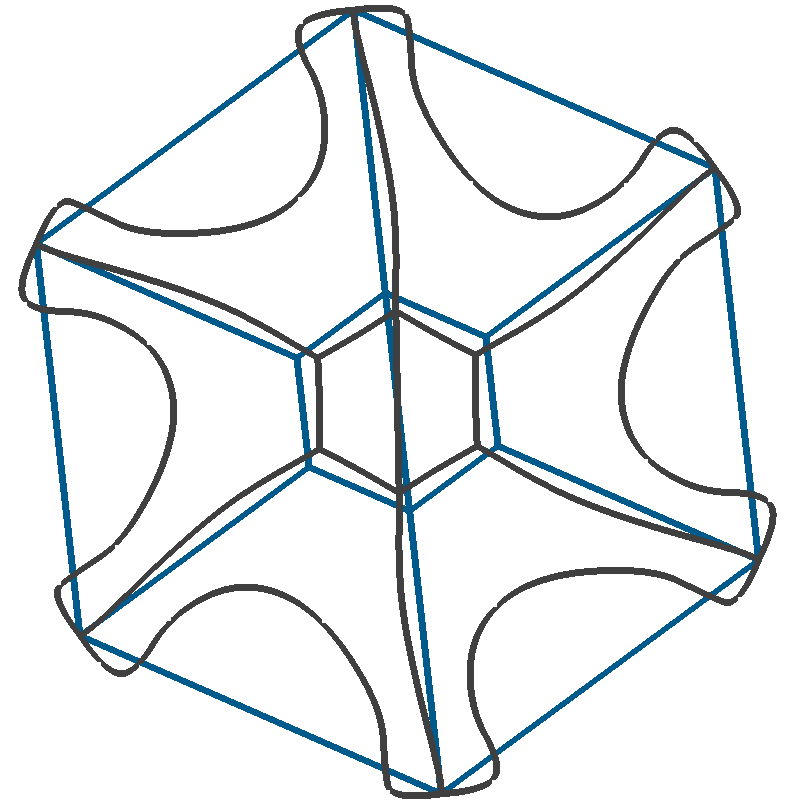}
	\caption{Multi-patch structures of the initial (left) and the deformed (right) domains.}\label{fig:femalecomp}
\end{figure}

Here we show that the mesh deformation approach is applicable to multi-patch problems as well. Consider the female rotor example depicted in Figure~\ref{fig:female}. The initial domain consists of 8 patches connected in a $C^0$-fashion. Each patch is formed by linear interpolation between its corner points, which corresponds to applying the $L^2$-simplification with the coarsest basis of degree 1. The mesh deformation is conducted by Newton's method with N-DIL with $N=5$ loading steps and Poisson's ratio of $0.48$. It is interesting to observe the way $C^0$-interfaces between the patches deform in an attempt to assume a more natural shape, see Fig.~\ref{fig:femalecomp}. 

\subsection{3D puzzle piece}
\begin{figure}[H]
	\centering
	\includegraphics[height=6.2cm]{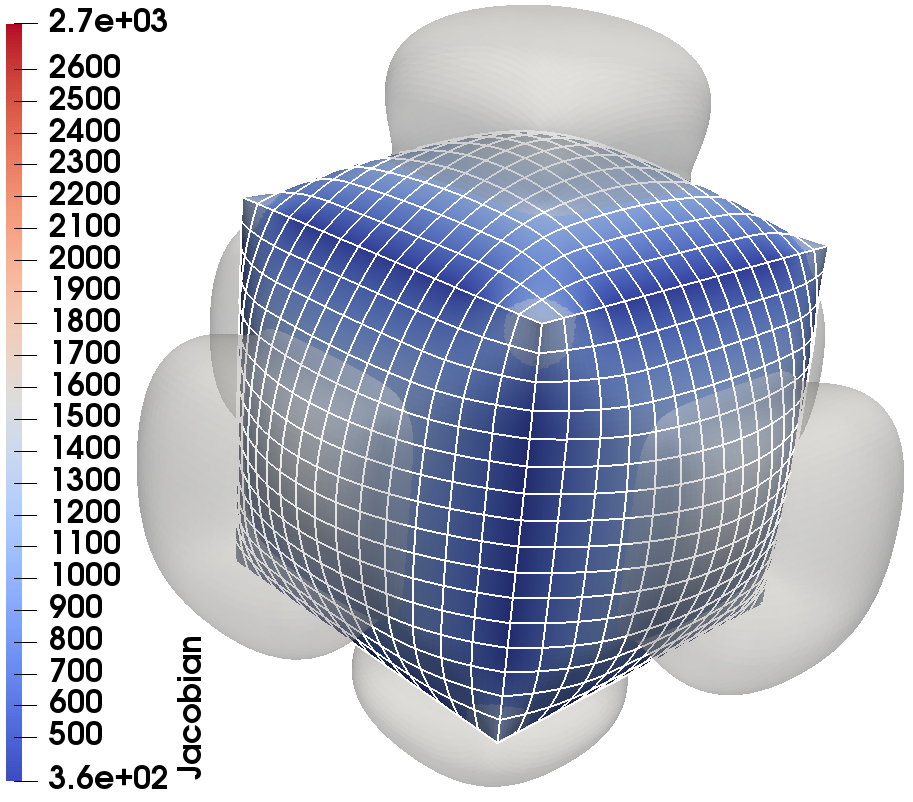}
	\includegraphics[height=6.2cm]{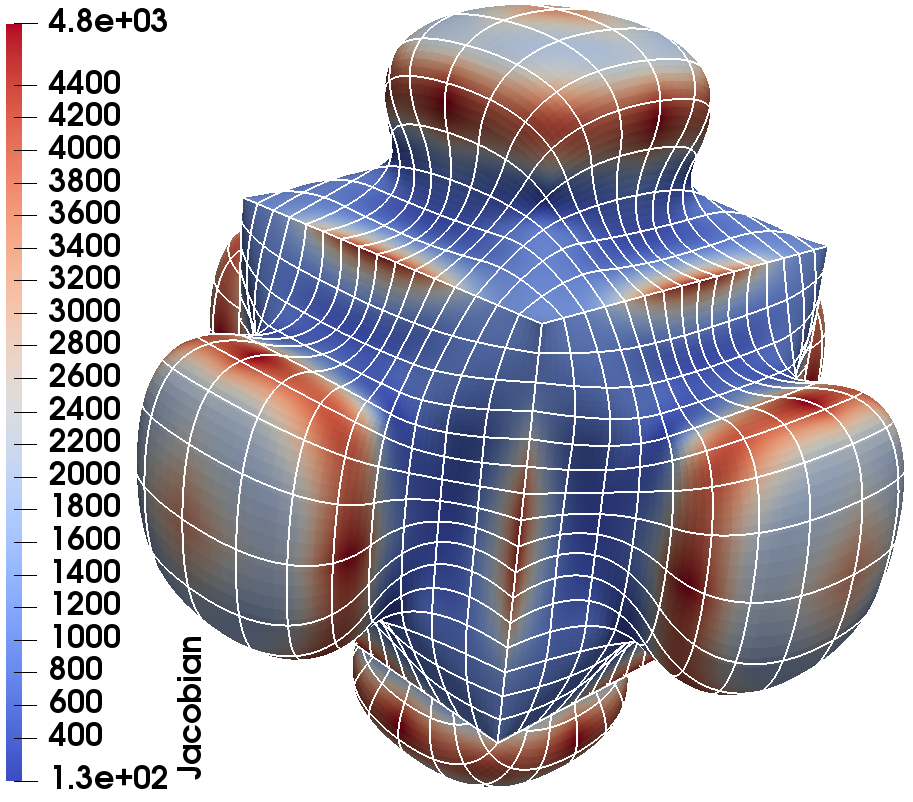}
	\caption{Mesh deformation approach for the 3D puzzle piece. Initial domain (left) and resulting parametrization (right).}\label{fig:puzzle3d}
\end{figure}

\begin{figure}[H]
	\centering
	\includegraphics[height=6.2cm]{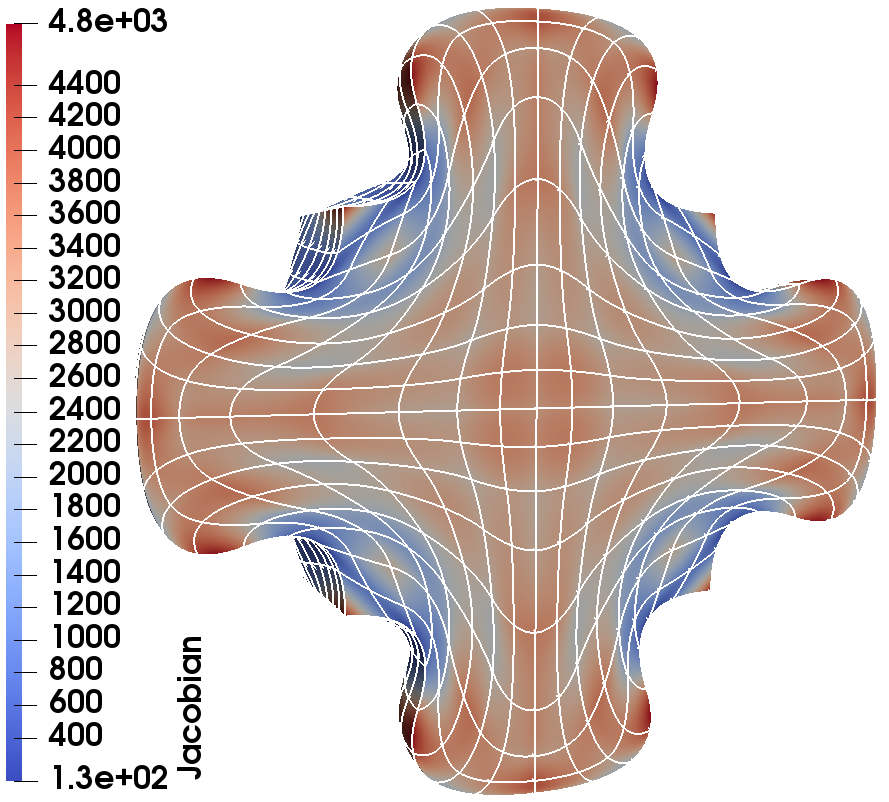}
	\caption{Cross-section of the 3D puzzle piece.}\label{fig:puzzle3dcross}
\end{figure}

Finally, we demonstrate that the mesh deformation approach is fully capable of dealing with 3D domains. Figures~\ref{fig:puzzle3d} and \ref{fig:puzzle3dcross} depict the result of applying it to a 3D puzzle piece example. The puzzle surface is simplified by the $L^2$-projection using the coarsest quadratic basis. The resulting initial domain is deformed using Newton's method with N-DIL with $N=10$ loading steps and Poisson's ratio equal to $0.46$.  

%% file: conclusion.tex
\section{Conclusion}
In this paper, we investigated the mesh deformation approach to the problem of domain parametrization and used it to construct tensor product B-spline parametrizations of high quality. We proposed a general technique to generate initial domains which can be applied to a wide range of examples. Furthermore, we described several efficient algorithms for computing an approximate solution to arising equations of nonlinear elasticity tuned specifically for this application. We demonstrated the performance of the mesh deformation approach on two 2D examples and compared it to the elliptic grid generation and area-orthogonality based optimization techniques. While being relatively computationally inexpensive, the proposed approach successfully produced bijective parametrizations which are superior with respect to uniformity of the corresponding mesh. Additionally, we showed that the mesh deformation approach is not restricted to a 2D single-patch case but can be applied to 3D and multi-patch problems.

Further research directions include development of an automatic procedure for the choice of an optimal initial domain. Moreover, the proposed approach may benefit from the nonhomogeneous distribution of material parameters in the elasticity model; potentially, a specialized material law can be developed. Lastly, a use of a nonzero right-hand side in the equations of nonlinear elasticity may offer more room for improvement.

\subsection*{Acknowledgement}
We thank Andreas Br\"{u}mmer and Matthias Utri for supplying us with the boundary profiles of the screw machine rotors. We also grateful to the team behind the G+Smo library for providing the access to an extensive collection of spline geometries from which the puzzle piece was selected. This research is supported by the German Research Council (DFG) under grant no. SI 756/5-1 (project YASON)  and by the German Federal Ministry of Education and Research (BMBF) under grant no. 05M16UKD  (project DYMARA).